\documentclass[11pt,reqno]{amsart}

\usepackage{fixltx2e} % Get the latest LaTeX enhancements
\usepackage{tgpagella} % Best font ever
\usepackage[utf8]{inputenc} % Customize to your encoding
\usepackage{amssymb, latexsym, stmaryrd, amsthm, dsfont, amsfonts, amsbsy, amsmath, mathrsfs} % Imprescindible
\usepackage{mathtools} % for some math typesetting tricks
\usepackage{bm, bbm} % for good bold math symbols
\usepackage{enumerate} % enhanced "enumerate" environments
\usepackage{verbatim} % for enhanced "comment" environment
\usepackage{url}   
\usepackage{lscape} % for \url command in bibliography
\usepackage{centernot}
\usepackage[dvipsnames]{xcolor}
\usepackage[colorinlistoftodos]{todonotes} % to add inline or margin comments

\usepackage[margin=1.1in]{geometry}

\usepackage{todonotes}

\usepackage{microtype,tikz,tikz-cd}  
\makeatletter % A correction needed for microtype:
\def\MT@register@subst@font{\MT@exp@one@n\MT@in@clist\font@name\MT@font@list
 \ifMT@inlist@\else\xdef\MT@font@list{\MT@font@list\font@name,}\fi}
\makeatother

\usepackage[pdftex,bookmarks,bookmarksnumbered,linktocpage,   %%  customize to your liking
         colorlinks,linkcolor=blue,citecolor=blue]{hyperref}

\newcommand{\bit}{\begin{itemize}}    % but see also \benbullet below
\newcommand{\eit}{\end{itemize}}
\newcommand{\ben}{\begin{enumerate}}
\newcommand{\een}{\end{enumerate}}
\newcommand{\benroman}{\ben[\normalfont (i)]}  % *
\let\eroman\een
\newcommand{\bde}{\begin{description}}
\newcommand{\ede}{\end{description}}
\theoremstyle{theorem}
\newtheorem{Theorem}{Theorem}[section]
\newtheorem{Theorem-n}{Theorem}

\newtheorem{Proposition}[Theorem]{Proposition}
\newtheorem{Modal Sahlqvist Theorem}[Theorem]{Modal Sahlqvist Theorem}
\newtheorem{Intuitionistic Sahlqvist Theorem}[Theorem]{Intuitionistic  Sahlqvist Theorem}
\newtheorem{Esakia Duality}[Theorem]{Esakia Duality}

\newtheorem{Main Lemma}[Theorem]{Main Lemma}

\newtheorem{Transfer Lemma}[Theorem]{Transfer Lemma}
\newtheorem{Abstract Sahlqvist Theorem}[Theorem]{Abstract Sahlqvist Theorem}

\newtheorem{Lemma}[Theorem]{Lemma}
\newtheorem{Corollary}[Theorem]{Corollary}
\newtheorem{Claim}[Theorem]{Claim}

\theoremstyle{definition}
\newtheorem{Definition}[Theorem]{Definition}
\newtheorem{exa}[Theorem]{Example}

\theoremstyle{remark}

\let\leq=\leqslant
\let\nleq=\nleqslant
\let\geq=\geqslant 

\usepackage[mathscr]{euscript}
 \let\mathscr\relax % just so we can load this and rsfs
\usepackage[scr]{rsfso}

\newcommand{\sub}{\subseteq}

\newcommand{\h}[1]{\mathsf{h}\left({#1}\right)}

\renewcommand{\int}{\mathsf{int}\,}

\bmdefine{\A}{A} 
\bmdefine{\C}{C}                                %  particular algebras
                              
\bmdefine{\2}{2}
\bmdefine{\B}{B}
\bmdefine{\D}{D}
\usepackage{scalerel}

\DeclareUnicodeCharacter{001B}{}
\subjclass[2020]{06D20, 06F30, 06E15, 54F05}

\keywords{prime spectrum, spectral space, Priestley space, Esakia space, distributive lattice, Heyting algebra, well-ordered tree, root system}

\begin{document}

\title[Trees and spectra of Heyting algebras]{Trees and spectra of Heyting algebras}

\author{Damiano Fornasiere and Tommaso Moraschini}

\address{Damiano Fornasiere and Tommaso Moraschini: Departament de Filosofia, Facultat de Filosofia, Universitat de Barcelona (UB), Carrer Montalegre, $6$, $08001$ Barcelona, Spain}

\email{damiano.fornasiere@gmail.com, tommaso.moraschini@ub.edu}

\date{\today}

\begin{abstract}
A poset is \emph{Esakia representable} when it is isomorphic to the prime spectrum of a Heyting algebra. Notably, every Esakia representable poset is also the spectrum of a commutative ring with unit. The problem of describing the Esakia representable posets was raised in 1985 and remains open to this day. We recall that a \emph{forest} is a disjoint union of trees and that a \emph{root system} is the order dual of a forest. It is shown that a root system is Esakia representable if and only if it satisfies a simple order theoretic condition, known as ``having enough gaps'', and each of its nonempty chains has an infimum.\ This strengthens Lewis's characterization of the root systems which are spectra of commutative rings with unit. While a similar characterization of arbitrary Esakia representable forests seems currently out of reach, we show that a \emph{well-ordered} forest is Esakia representable if and only if it has enough gaps and each of its nonempty chains has a supremum.
\end{abstract}

\maketitle

%\tomnote{Why is Speed reference precise and Joyal isn't? Why ``see also...'' in Speed?}

\section{Introduction}

%\tomnote{Refer to recentg results by Wherung?}

The \emph{prime spectrum} of a commutative ring with unit is the poset of its prime ideals.\ The \emph{representation problem}, raised by Kaplansky in \cite[pp.\ 5--7]{Ka74}, asks for a characterization of the posets isomorphic to the prime spectra of commutative rings with unit.\ A similar problem was raised by Gr\"atzer in \cite[Problem 34, p.\ 156]{Gr71} in the context of order theory. We recall that the \emph{prime spectrum} of a bounded distributive lattice is the poset of its prime filters. Gr\"atzer's problem asks for a characterization of the posets isomorphic to the prime spectra of bounded distributive lattices. 

Notably, the two problems coincide because commutative rings with unit and bounded distributive lattices have the same prime spectra (see,
e.g., \cite[Thm.\ 1.1]{Priestley94sp}).\ More precisely, Hochster showed that the prime spectra
of commutative rings with unit endowed with the Zariski topology are precisely the spectral spaces \cite{Hochster69} (see also \cite{specDST85}) and Stone did the same for the prime spectra of bounded distributive lattices  \cite{St38a}. Because of this, we say that a poset is \emph{representable} when it is isomorphic to the prime spectrum of a commutative ring with unit (equiv.\ of a bounded distributive lattice).\ In this parlance, the representation problem asks for a characterization of the representable posets.

Some conditions equivalent to the representability of a poset are known. For instance,  Joyal \cite{JoyalFeb71} and Speed \cite[Thm.\ p.\ 85]{Speed72}  showed that a poset is representable if and only if it is profinite. Moreover, in view of Priestley duality \cite{Pr70,Pr72}, a poset is representable precisely when it can be endowed with a topology that turns it into a Priestley space. However, these characterizations provide little information on the inner structure of representable posets, which is why the representation problem remains elusive to this day. 

One of the main positive results on the inner structure of representable posets is due to Lewis \protect{\cite[Thm.\ 3.1]{Lwevis73}}. We recall that a poset is a \emph{tree} when it is rooted and its principal downsets are chains and that it is a \emph{root system} when it is a disjoint union of order duals of trees. Lewis showed that a root system $X$ is representable if and only if each of its nonempty chains has an infimum and $X$ has \emph{enough gaps}, where the latter means that if $x < y$, there exist $z, v \in [x, y]$ such that $z$ is an immediate predecessor of $v$. Since the class of representable posets is closed under the formation of order duals, we obtain that a \emph{forest} (i.e., a disjoint union of trees) is representable if and only if it has enough gaps and each of its nonempty chains has a supremum.

In this paper, we focus on the representation problem for a prominent class of bounded distributive lattices, namely, \emph{Heyting algebras}.\ These are the bounded distributive lattices in which the meet operation has an adjoint, sometimes called implication.\ Heyting algebras arise in different areas of mathematics, including:

\benroman
\item \textsf{topology}: the lattice of open sets of any topological space is a Heyting algebra;
\item \textsf{domain theory}: each continuous distributive lattice is a Heyting algebra;
\item \textsf{topos theory}: the subobject classifier of any topos is a Heyting algebra;
\item \textsf{algebra}: any distributive algebraic lattice is a Heyting algebra;
\item \textsf{constructivism and logic}: the algebraic models of intuitionistic logic are Heyting algebras;
\item \textsf{order theory}: Heyting algebras are the most common generalization of Boolean algebras.
\eroman

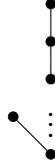
\begin{figure}
\[
	\begin{tikzpicture}
		\tikzstyle{point} = [shape=circle, thin, draw=black, fill=black, scale=0.35]
		
		\node (0) at (0,0)[point] {};
		\node (1) at (0,-0.5)[point] {};
		\node (2) at (0,-1)[point] {};
		\node (d) at (0,-1.5)[] {$\vdots$};
		
		\node (omega) at (0,-2)[point] {};
		\node (omega1) at (-0.5,-1.5)[point]{};
		
		\draw (0)--(1)--(2);
		\draw (omega)--(omega1);
	\end{tikzpicture}
	\]
\caption{A representable tree that is not Esakia representable.}\label{Figure 1}
\end{figure}

In 1985, Esakia raised the problem of describing the posets isomorphic to the prime spectra of Heyting algebras \cite[Appendix A.5]{Esakia-book85}. Accordingly, we say that a poset is \emph{Esakia representable} when it is isomorphic to the prime spectrum of a Heyting algebra. While every Esakia representable poset is representable in the traditional sense (because every Heyting algebra is a bounded distributive lattice), the converse does not hold in general: for instance, the poset depicted in Figure \ref{Figure 1}  is representable, but not Esakia representable (see \cite[Example 5.6]{BeMo09}).

After four decades, the problem of describing the Esakia representable posets remains open. In addition, this problem cannot be reduced to the one of describing the representable posets because no concrete way of isolating the Esakia representable posets from the class of all the representable posets is known. Nonetheless, some progress has been made and, recently, a characterization of the Esakia representable root systems whose maximal chains are either finite or of order type dual to $\omega+1$ was obtained (see \cite[Cor.\ 6.20]{GNTSAiM11} and its proof).

In this paper, we extend this result by providing a description of all the Esakia representable root systems. More precisely, we show that a root system is Esakia representable if and only if has enough gaps and each of its nonempty chains has an infimum (Theorem \ref{Thm:main-root}). As a corollary, we obtain Lewis' classical description of the representable root systems. We recall that the Heyting algebras whose prime spectrum is a root system have been called \textit{G\"odel algebras} \cite{Ha98} (see also \cite[Thm.\ 2.4]{Ho69}). Therefore, our result takes the form of a characterization of the prime spectra of G\"odel algebras.

Contrarily to the case of arbitrary representable posets, the class of Esakia representable posets is not closed under order duals. In particular, the tree depicted in Figure \ref{Figure 1} is not Esakia representable, although its order dual is Esakia representable because it has enough gaps and its nonempty chains have infima.\ Notice that the tree in Figure \ref{Figure 1} contains an infinite descending chain. We will show that Lewis' description of the representable forests can be extended to Esakia representable forests by prohibiting the presence of such chains. More precisely, a forest is said to be \emph{well-ordered} when it has no infinite descending chain. We show that a well-ordered forest is Esakia representable if and only if each of its nonempty chains has a supremum (Theorem \ref{Thm:main-forest}).

At the heart of our proof stands a novel compactness argument which combines intuitions from combinatorics, algebra, and topology and highlights the higher complexity of Esakia representable forests, as opposed to arbitrary representable forests. It remains an  open problem to give a full characterization of arbitrary (i.e., not necessarily well-ordered) Esakia representable forests.

\section{Representable posets}

Among bounded distributive lattices, a special role is played by Heyting algebras \cite{BaDw74,Esakia-book85,RaSi63}. We recall that a bounded distributive lattice $A$ is said to be a \textit{Heyting algebra} when it can be expanded with a binary operation $\to$ such that for every $a, b, c \in A$,
\[
a \land b \leq c \, \, \Longleftrightarrow \, \, a \leq b \to c.
\]
In this case, this expansion is unique and $b \to c = \max \{ a \in A : a \land b \leq c \}$.

In view of \textit{Priestley} and \textit{Esakia dualities} \cite{Es74,Esakia-book85,Pr70,Pr72}, the problem of describing the spectra of bounded distributive lattices and Heyting algebras can be phrased in purely topological terms, as we proceed to illustrate.  Given a poset $\langle X; \leq \rangle$ and $Y \subseteq X$, let
\[
{\uparrow} Y \coloneqq \{x\in X : \exists y \in Y \text{ s.t. } y\leq x\} \, \, \text{ and } \, \, {\downarrow} Y \coloneqq \{x\in X : \exists y \in Y \text{ s.t. } x\leq y\}.
\]
The set $Y$ is said to be an {\em upset} (resp.\ {\em downset}) if $Y={\uparrow}Y$ (resp.\ $Y={\downarrow}Y$). When $Y = \{ x \}$, we will write ${\uparrow} x$ and ${\downarrow} x$ instead of ${\uparrow} \{ x \}$ and ${\downarrow} \{ x\}$. We denote the set of clopen upsets of an ordered topological space $X$ by $\mathsf{ClUp}\left(X\right)$.

\begin{Definition}
An ordered topological space $X = \langle X; \leq, \tau \rangle$ is a \textit{Priestley space} when it is compact and satisfies the \textit{Priestley separation axiom}: for every $x, y \in X$,
\[
\text{if }x \nleq y, \text{ there exists }U \in \mathsf{ClUp}\left(X\right) \text{ such that }x \in U \text{ and }y \notin U.
\]
If, in addition, ${\downarrow}U$ is open for every open set $U \subseteq X$, then $X$ is said to be an \textit{Esakia space}.
\end{Definition}

Let $A$ be a bounded distributive lattice.\ A set $F \subseteq A$ is a \emph{prime filter} of $A$ when it is a nonempty proper upset such that for every $a, b \in A$,
\[
(a, b \in F \Longrightarrow a \land b \in F) \, \, \text{ and } \, \, (a \lor b \in F \Longrightarrow a \in F \text{ or }b \in F).
\]
The \emph{prime spectrum} of $A$ is the poset $\langle \mathsf{Pr}\left(A\right); \subseteq \rangle$, where $\mathsf{Pr}\left(A\right)$ is the set of prime filters of $A$.

Now, for each $a \in A$ let
\[
\gamma_A\left(a\right) \coloneqq \{ F \in \mathsf{Pr}\left(A\right) : a \in F \}.
\]
Then the triple $A_+ \coloneqq \langle \mathsf{Pr}\left(A\right); \subseteq, \tau \rangle$, where  $\tau$ is the topology on $\mathsf{Pr}\left(A\right)$ generated by the subbase
\[
\{ \gamma_A\left(a\right) : a \in A \} \cup \{ \gamma_A\left(a\right)^c : a \in A \},
\]
is a Priestley space.

If, moreover, $A$ is a Heyting algebra, then $A_+$ is an Esakia space. On the other hand, given a Priestley space $X$, the structure $X^+ \coloneqq \langle \mathsf{ClUp}\left(X\right); \cap, \cup, \emptyset, X \rangle$ is a bounded distributive lattice. If, in addition, $X$ is an Esakia space, then $X^+$ is a Heyting algebra in which the operation $\to$ is defined as
\[
U \to V \coloneqq \{ x \in X : U \cap {\uparrow}x \subseteq V \}.
\]
It is a consequence of Priestley and Esakia dualities that these transformations are one inverse to the other, in the sense that
\begin{equation}\label{Eq : the two natural isomorphisms}
A \cong \left(A_+\right)^+ \, \, \text{ and } \, \, X \cong \left(X^+\right)_+.
\end{equation}

Recall that a poset is \emph{representable} (resp. \emph{Esakia representable}) when it is isomorphic to the prime spectrum of a bounded distributive lattice (resp.\ Heyting algebra). The next observation is a consequence of the isomorphisms in condition (\ref{Eq : the two natural isomorphisms}).

\begin{Theorem}\label{Thm:basic}
The following conditions hold:
\benroman
\item\label{basic:item:1} A poset is representable iff it can be endowed with a topology that turns it into a Priestley space;
\item\label{basic:item:2} A poset is Esakia representable iff it can be endowed with a topology that turns it into an Esakia space.
\eroman
\end{Theorem}

While the structure of (Esakia) representable posets remains largely unknown, they need to satisfy a number of nontrivial properties. Given a poset $\langle X; \leq \rangle$ and $x, y \in X$, we say that $x$ is an \textit{immediate predecessor} of $y$ when $x < y$ and there exists no $z \in X$ such that $x < z < y$. We write $x \prec y$ to indicate that this is the case.

\begin{Definition}
A poset $X = \langle X; \leq \rangle$ is said to
\benroman
\item  have \emph{enough gaps} when for every $x, y \in X$ such that $x < y$, there exist $x' \geq x$ and $y' \leq y$ such that $x' \prec y'$;
\item be \textit{Dedekind complete} when every nonempty chain in $X$ has a supremum and an infimum.
\eroman
\end{Definition}

A subset $U$ of a poset $X$ is \textit{order open} when it belongs to the least family $\mathcal{O}$ of subsets of $X$ such that:
\benroman
\item $\{ x \}^c \in \mathcal{O}$ for every $x \in X$;
\item\label{order-closed:2} if $U \in \mathcal{O}$, then $\left({\uparrow}\left(U^c\right)\right)^c, \left({\downarrow} \left(U^c\right)\right)^c \in \mathcal{O}$;
\item\label{order-closed:3} $\mathcal{O}$ is closed under finite intersections and arbitrary unions.
\eroman

\begin{Definition}
A poset $X$ is said to be \emph{order compact} when for every  family $\{ U_i : i \in I \}$ of order open sets,
\[
\text{if }\bigcup_{i \in I} U_{i} = X, \text{ there exists a finite }J \subseteq I \text{ such that }\bigcup_{j \in J} U_{j} = X.
\]
\end{Definition}

\begin{Proposition}\label{Prop:properties}
Representable posets have enough gaps and are both Dedekind complete and order compact.
\end{Proposition}

\begin{proof}
For the fact that representable posets have enough gaps and are Dedekind complete, see \cite[pp.\ 5--7]{Ka74}. On the other hand, every representable poset $X$ is order compact because the order open sets are open in any topology that turns $X$ into a Priestley space and Priestley spaces are compact (a slightly weaker statement can be found in \cite[p.\ 822, condition (H)]{LewOhm76}.
\end{proof}

The converse of Proposition \ref{Prop:properties} does not hold, however, as shown in  \cite[Example 2.1]{LewOhm76}. We will also rely on the following observation.

\begin{Proposition}\label{Prop:closure}
The following conditions hold:
\benroman
\item\label{Prop:closure:item:1} The class of representable posets is closed under disjoint unions and order duals;
\item\label{Prop:closure:item:2} The class of Esakia representable posets is closed under disjoint unions.
\eroman
\end{Proposition}
\begin{proof}
Condition (\ref{Prop:closure:item:2}) is Proposition 5.1(6) in the Appendix of \cite{Esakia-book85}. Therefore, we turn to prove Condition (\ref{Prop:closure:item:1}). The fact that the class of representable posets is closed under order duals follows immediately from Theorem \ref{Thm:basic}(\ref{basic:item:1}) and the fact that if $\langle X; \leq, \tau \rangle$ is a Priestley space, so is $\langle X; \geq, \tau \rangle$. On the other hand, closure under disjoint unions holds by \cite[Thm.\ 4.1]{LewOhm76}.
\end{proof}	

Notice that the class of Esakia representable posets is not closed under order duals because the poset in Figure \ref{Figure 1} is not Esakia representable \cite[Example 5.6]{BeMo09}, although its order dual is \cite[Cor.\ 6.20]{GNTSAiM11}.

Given a pair of sets $X$ and $Y$, we will write $X \sub_\omega Y$ to indicate that $X$ is a finite subset of $Y$. We will rely on the following easy observation.

\begin{Lemma}\label{Lem : typical order open sets}
Let $X$ be a poset and $Y, Z \sub_\omega X$. Then $\left({\uparrow}Y \cap {\downarrow}Z\right)^c$ is an order open set of $X$.
\end{Lemma}

\begin{proof}
We will show that $Y^c$ and $Z^c$ are order open. By symmetry it suffices to prove that  $Y^c$ is order open. If $Y = \emptyset$, then $Y^c = X$ is the intersection of the empty family. As the family of order open sets is closed under finite (possibly empty) intersections, we are done. Then we consider the case where $Y \ne \emptyset$. Consider an enumeration $Y= \{ y_1, \dots, y_n \}$. Since the sets $\{ y_1\}^c, \dots, \{ y_n \}^c$ are order open, so is their intersection $Y^c = \{ y_1\}^c \cap \dots\cap \{ y_n \}^c$. Hence, $Y^c$ and $Z^c$ are order open sets as desired. As a consequence, $\left({\uparrow}Y\right)^c$ and $\left({\downarrow}Z\right)^c$ are order open too and so is their union $\left({\uparrow}Y\right)^c  \cup \left({\downarrow}Z\right)^c$. Since $\left({\uparrow}Y\right)^c  \cup \left({\downarrow}Z\right)^c = \left({\uparrow}Y \cap {\downarrow}Z\right)^c$, we are done.
\end{proof}

Throughout the paper, we denote the class of all ordinals by $\mathsf{Ord}$. Furthermore, given a poset $ \langle X; \leq \rangle$ and a set $Y \subseteq X$, we denote the sets of maximal and minimal elements of the subposet $\langle Y; \leq \rangle$ by $\max Y$ and $\min Y$, respectively. Furthermore, when they exist, we let $\sup Y$ and $\inf Y$ be the supremum and infimum of $Y$, respectively.

\section{Esakia representable root systems}
	
\begin{Definition}
A poset is said to be 
\benroman
\item a \textit{tree} when it is rooted and its principal downsets are chains;
\item a \textit{forest} when it is isomorphic to the disjoint union of a family of trees;
\item a \textit{root system} when it is the order dual of a forest.
\eroman
\end{Definition}	
	
One of the main positive results on the representation problem is the next theorem of Lewis.

\begin{Theorem}[\protect{\cite[Thm.\ 3.1]{Lwevis73}}]\label{Thm:Lewis-root}
A root system is representable iff it has enough gaps and each of its nonempty chains has an infimum.
\end{Theorem}

In this section, we strengthen this result by showing that it still holds in the context of Esakia representable posets. To this end, we recall that a Heyting algebra is a \textit{G\"odel algebra} \cite{Ha98} when it validates the equation
\[
\left(x \to y\right) \lor \left(y \to x\right) \thickapprox 1
\]
or, equivalently, it is isomorphic to a subdirect product of chains \cite[Thm.\ 1.2]{Ho69}. From a logical standpoint, the importance of G\"odel algebras comes from the fact that they algebraize the \textit{G\"odel-Dummett logic} \cite{Dm59} in the sense of \cite{BP89} (see, e.g., \cite{ChZa97}). Notably, G\"odel algebras can be characterized in term of the shape of their spectra.

\begin{Theorem}[\protect{\cite[Thm.\ 2.4]{Ho69}}]\label{Thm:Horn}
A Heyting algebra is a G\"odel algebra iff its prime spectrum is a root system.
\end{Theorem}

From Theorems \ref{Thm:basic}(\ref{basic:item:2}) and \ref{Thm:Horn} we deduce the following.

\begin{Corollary}
A poset is isomorphic to the prime spectrum of a G\"odel algebra iff it is an Esakia representable root system.
\end{Corollary}

The aim of this section is to establish the following description of the Esakia representable root systems (equiv.\ of the prime spectra of G\"odel algebras).

\begin{Theorem}\label{Thm:main-root}
A root system is Esakia representable iff it has enough gaps and each of its nonempty chains has an infimum.
\end{Theorem}

We remark that Theorem \ref{Thm:Lewis-root} is an immediate consequence of Theorem \ref{Thm:main-root}. More precisely, the implication from left to right in Theorem \ref{Thm:Lewis-root} holds by Proposition \ref{Prop:properties}, while the other implication holds by Theorem \ref{Thm:main-root} and the fact that every Esakia  representable poset is  representable. Furthermore, a weaker version of Theorem \ref{Thm:main-root}, stating that the result holds for the root systems whose maximal chains are either finite or of order type dual to $\omega+1$, can be deduced from \cite[Cor.\ 6.20]{GNTSAiM11}.

\begin{proof}[Proof of Theorem \ref{Thm:main-root}]
In view of Proposition \ref{Prop:properties}, it suffices to prove the implication from right to left. To this end, it will be enough to show that the following condition holds for every poset $X$ whose order dual is a tree:
\begin{align}\label{Eq:root-system-1}
    \begin{split}
    \text{if }X&\text{ has enough gaps and each of its nonempty chains has an infimum,}\\
   & \text{ then $X$ is Esakia representable.}
\end{split}
\end{align}
For suppose that condition (\ref{Eq:root-system-1}) holds for the order duals of trees and consider a root system $X$ with enough gaps and in which each nonempty chain has an infimum. Since $X$ is a root system, it is the disjoint union of a family of posets $\{ X_i : i \in I \}$ whose order duals are trees. Furthermore, each $X_i$ has enough gaps as well as infima of nonempty chains. Therefore, each $X_i$ is Esakia representable by condition (\ref{Eq:root-system-1}). Hence, the disjoint union $X$ is also Esakia representable by Proposition \ref{Prop:closure}(\ref{Prop:closure:item:2}) as desired.

Therefore, we turn to prove condition (\ref{Eq:root-system-1}). Consider a poset $X = \langle X; \leq \rangle$ with enough gaps, in which every nonempty chain has an infimum, and whose order dual is a tree. Let then $\tau$ be the topology on $X$ generated by the subbase
\[
\mathcal{S} \coloneqq \{ {\downarrow}x : \exists y \in X \text{ s.t. }x \prec y
 \} \cup \{ \left({\downarrow}x\right)^c : \exists y \in X \text{ s.t. }x \prec y \}.
\]
We will show that $X = \langle X; \leq, \tau \rangle$ is an Esakia space. The proof proceeds through a series of claims.

\begin{Claim}\label{Claim:compact1}
The topological space $X$ is compact.
\end{Claim}

\begin{proof}[Proof of the Claim.]
Suppose the contrary, with a view to contradiction. By Alexander's subbase theorem there exists an open cover $\mathcal{C} \subseteq \mathcal{S}$ of $X$ without any finite subcover.  To this end, we will define recursively a sequence $\{ x_\alpha : \alpha \in \mathsf{Ord} \}$ of elements of $X$ such that for every ordinal $\alpha$,
\benroman
\item\label{1c} $\left({\downarrow}x_\alpha\right)^c \in \mathcal{C}$;
\item\label{2c} $x_\beta < x_\gamma$ for every $\gamma < \beta \leq \alpha$.
\eroman 
Clearly, the validity of condition (\ref{2c}) for every ordinal $\alpha$ implies that $X$ is a proper class, which is the desired contradiction.

Consider an ordinal $\alpha$ and suppose that we already defined a sequence $\{ x_\beta : \beta < \alpha\}$ of elements of $X$ such that
\benroman
\renewcommand{\labelenumi}{(\theenumi)}
\renewcommand{\theenumi}{L\arabic{enumi}}
\item\label{item : limit : L : 1} $\left({\downarrow}x_\beta\right)^c \in \mathcal{C}$ for each $\beta < \alpha$;
\item\label{item : limit : L : 2}  $x_\beta < x_\gamma$ for every $\gamma < \beta < \alpha$.
\eroman
We will prove that the set $Y \coloneqq \{ x_\beta : \beta < \alpha\}$ has an infimum in $X$. If $Y = \emptyset$, then $\inf Y$ is the maximum of $X$, which exists because $X$ is the order dual of a tree. The we consider the case where $Y \ne \emptyset$. In view of condition (\ref{item : limit : L : 2}), the set $Y$ is a chain. As nonempty chains have infima by assumption, we conclude that $\inf Y$ exists.

Since $\mathcal{C}$ covers $X$, there exists $U \in \mathcal{C}$ such that $\inf Y \in U$. Furthermore, as $\mathcal{C} \subseteq \mathcal{S}$, there also exists $z \in X$ such that 
\benroman
\renewcommand{\labelenumi}{(\theenumi)}
\renewcommand{\theenumi}{C\arabic{enumi}}
\item\label{item : z is top element in Priestley case} $z$ has an immediate successor;
\item\label{item : 2 z is top element in Priestley case} either $U = {\downarrow}z$ or $U = \left({\downarrow}z\right)^c$.
\eroman

We will show that the case where $U = {\downarrow}z$ never happens. Suppose the contrary, with a view to contradiction. We have two cases: either $\inf Y \in Y$ or $\inf Y \notin Y$. First, suppose that $\inf Y \in Y$. Since $\inf Y \in U = {\downarrow}z$, we have $X = U \cup ({\downarrow}\inf Y)^c$. As $U \in \mathcal{C}$ and $\mathcal{C}$ lacks a finite subcover by assumption, this yields $({\downarrow}\inf Y)^c \notin \mathcal{C}$. On the other hand, from $\inf Y \in Y$ and condition (\ref{item : limit : L : 1})  it follows that $({\downarrow}\inf Y)^c \in \mathcal{C}$, a contradiction. Then we consider the case where $\inf Y \notin Y$.
Together with the fact that ${\uparrow}\inf Y$ is a chain (because the order dual of $X$ is a tree), this implies that $\inf Y$ does not have immediate successors.  By condition (\ref{item : z is top element in Priestley case}) we obtain $\inf Y \ne z$. Therefore, from $\inf Y \in U = {\downarrow}z$ it follows that $\inf Y < z$. As ${\uparrow}\inf Y$ is a chain and $Y = \{ x_\beta : \beta < \alpha\}$, there exists $\beta < \alpha$ such that $x_\beta < z$. By condition (\ref{item : limit : L : 1}) we have $\left({\downarrow}x_\beta\right)^c \in \mathcal{C}$ which, together with $x_\beta < z$ and ${\downarrow}z = U \in \mathcal{C}$, implies that $\{\left({\downarrow}x_\beta\right)^c, {\downarrow}z\}$ is a finite subcover of $\mathcal{C}$, a contradiction. Therefore, we conclude that $U \ne {\downarrow}z$ as desired. By condition (\ref{item : 2 z is top element in Priestley case}) this means that
 $U = \left({\downarrow}z\right)^c$.
 
  We will prove that $z < x_\beta$ for every $\beta < \alpha$. Suppose, on the contrary, that there exists $\beta < \alpha$ such that $z = x_{\beta}$ or $z \nleq x_\beta$. From $\inf Y \in U = \left({\downarrow}z\right)^c$ it follows that $\inf Y \nleq z$. Since $x_\beta \in Y$, this yields $x_\beta \nleq z$. Together with the assumption that either $z = x_{\beta}$ or $z \nleq x_\beta$, this implies $z \nleq x_\beta$. Consequently, $x_\beta$ and $z$ are incomparable. As the order dual of $X$ is a tree, this guarantees that ${\downarrow}x_\beta \cap {\downarrow}z = \emptyset$. Hence,
\[
\left({\downarrow}x_\beta\right)^c \cup \left({\downarrow}z\right)^c = \left({\downarrow}x_\beta \cap {\downarrow}z\right)^c = \emptyset^c = X.
\]
Since $\left({\downarrow}z\right)^c = U \in \mathcal{C}$ and $\left({\downarrow}x_\beta\right)^c \in \mathcal{C}$ (the latter by condition (\ref{item : limit : L : 1})), we obtain that $\{ \left({\downarrow}x_\beta\right)^c, \left({\downarrow}z\right)^c\}$ is a finite subcover of $\mathcal{C}$, a contradiction.\ Hence, we conclude that $z < x_\beta$ for every $\beta < \alpha$.  Thus, letting $x_\alpha \coloneqq z$, we obtain  $x_\alpha < x_\beta$ for every $\beta < \alpha$. Since $\left({\downarrow}x_\alpha\right)^c = \left({\downarrow}z\right)^c = U \in \mathcal{C}$,  the elements in the sequence $\{ x_\beta : \beta \leq \alpha \}$ satisfy conditions (\ref{1c}) and (\ref{2c}) as desired.

This completes the recursive definition of the sequence $\{ x_\alpha: \alpha \in \mathsf{Ord} \}$ and produces the desired contradiction.
\end{proof}

\begin{Claim}\label{Claim:Priestley1}
The ordered topological space $X$ satisfies Priestley separation axiom.
\end{Claim}

\begin{proof}[Proof of the Claim.]
Consider $x, y \in X$ such that $x \nleq y$. If $y$ has an immediate successor, we have  ${\downarrow}y, \left({\downarrow}y\right)^c \in \mathcal{S}$ by the definition of $\mathcal{S}$. In this case, $\left({\downarrow}y\right)^c$ is a clopen upset containing $x$ and missing $y$ as desired. Then we consider the case where $y$ does not have immediate successors. Notice that $y$ is not the maximum of $X$, otherwise we would have $x \leq y$, which is false. Therefore, ${\uparrow}y \smallsetminus \{ y \} \ne \emptyset$. Furthermore, ${\uparrow}y \smallsetminus \{ y \}$ is a chain because the order dual of $X$ is a tree.\ Now, since ${\uparrow}y \smallsetminus \{ y \}$ is a nonempty chain, it has an infimum by assumption. As $y$ lacks immediate successors, this infimum must be $y$ itself. As a consequence, from $x \nleq y$ it follows that there exists $z > y$ such that $x \nleq z$. As $X$ has enough gaps, there exists also an element $y^+ \in X$ with an immediate successor and such that $y \leq y^+ < z$. Consequently, ${\downarrow}y^+, \left({\downarrow}y^+\right)^c \in \mathcal{S}$ by the definition of $\mathcal{S}$. Furthermore, $x \nleq y^+$ because $x \nleq z$ and $y^+ \leq z$. Thus, $\left({\downarrow}y^+\right)^c$ is a clopen upset containing $x$ and missing $y$.
\end{proof}

From Claims \ref{Claim:compact1} and \ref{Claim:Priestley1} it follows that $X$ is a Priestley space. In order to prove that it is also an Esakia space, we need to show that the downset of every open set is also open. To this end, let $\mathcal{B}$ be the base for the topology of $X$ consisting of all the finite intersections of the elements of the subbase $\mathcal{S}$. As every open set $U$ is the union of a family $\{ U_i : i \in I \} \subseteq \mathcal{B}$ and
\[
{\downarrow}U = \bigcup_{i \in I}{\downarrow}U_i,
\]
it will be enough to prove that the downset of every element of $\mathcal{B}$ is open.

Consider $U_1, \dots, U_n \in \mathcal{S}$. We need to show that ${\downarrow}\left(U_1 \cap \dots \cap U_n\right)$ is open. We may assume that $U_1 \cap \dots \cap U_n \ne \emptyset$, otherwise ${\downarrow}\left(U_1 \cap \dots \cap U_n\right) = \emptyset$ and we are done. By the definition of $\mathcal{S}$ for every $m \leq n$ there exists $x_m \in X$ such that either $U_m = {\downarrow}x_m$ or $U_m = \left({\downarrow}x_m\right)^c$. Let $Y \coloneqq \{ x_m : U_m = {\downarrow}x_m \}$ and let $Y^c$ be the complement of $Y$ relative to $\{ x_m : m \leq n \}$. Observe that
\begin{equation}\label{Eq:base-finite-intersection}
U_1\cap \dots \cap U_n = \bigcap_{x_m \in Y} {\downarrow}x_m \cap \bigcap_{y_m \in Y^c}\left({\downarrow}x_m\right)^c = \bigcap_{x_m \in Y} {\downarrow}x_m \cap \left({\downarrow}\left(Y^c\right)\right)^c.
\end{equation} 

We have two cases: either $Y = \emptyset$ or $Y\ne \emptyset$. First, suppose that $Y = \emptyset$. In view of the above equalities,  we have $U_1\cap \dots \cap U_n = \left({\downarrow}\left(Y^c\right)\right)^c$. As $U_1 \cap \dots \cap U_n \ne \emptyset$ by assumption, the upset   $\left({\downarrow}\left(Y^c\right)\right)^c$ is nonempty and, therefore, contains the maximum $\top$ of $X$. Consequently, $\top \in U_1\cap \dots \cap U_n$ and, therefore, ${\downarrow}(U_1 \cap \dots \cap U_n) = X$ is an open set. 

Then we consider the case where $Y \ne \emptyset$. We will prove that $Y$ is a chain. For if $Y$ contained two incomparable elements $x_{k}$ and $x_{m}$, we would have
\[
U_1\cap \dots \cap U_n \subseteq U_{k} \cap U_{m} = {\downarrow}x_{k} \cap {\downarrow}x_{m} = \emptyset,
\]
where the last equality follows from the assumption that $x_k$ and $x_m$ are incomparable and the order dual of $X$ is a tree. But this contradicts the assumption that $U_1\cap \dots \cap U_n \ne \emptyset$.

Now, since $Y$ is a finite nonempty chain, it has a minimum $y$. Consequently, condition (\ref{Eq:base-finite-intersection}) can be simplified as follows:
\begin{equation}\label{Eq:base-finite-intersection2}
U_1\cap \dots \cap U_n = {\downarrow}y \cap \left({\downarrow}Y^c\right)^c.
\end{equation}

We will prove that $y \in U_1\cap \dots \cap U_n$. In view of Condition (\ref{Eq:base-finite-intersection2}), it suffices to show that $y \in \left({\downarrow}Y^c\right)^c$. Suppose the contrary, with a view to contradiction. Then there exists $x_m \in Y^c$ such that $y \leq x_m$.  Consequently, ${\downarrow}y \cap \left({\downarrow}x_m\right)^c = \emptyset$. Together with condition (\ref{Eq:base-finite-intersection2}) and $x_m \in Y^c$, this implies $U_1\cap \dots \cap U_n = \emptyset$, a contradiction. Hence, we conclude that $y \in U_1\cap \dots \cap U_n$.

As a consequence, we obtain that ${\downarrow}y \subseteq {\downarrow}\left(U_1\cap \dots \cap U_n\right)$. Since the reverse inclusion holds by condition (\ref{Eq:base-finite-intersection2}), we conclude that ${\downarrow}\left(U_1\cap \dots \cap U_n\right) = {\downarrow}y$. From $y \in Y$ and the definition of $Y$ it follows that ${\downarrow}y = U_m$ for some $m \leq n$. Therefore, ${\downarrow}\left(U_1\cap \dots \cap U_n\right) = {\downarrow}y = U_m$. As $U_m \in \mathcal{S}$, we conclude that ${\downarrow}\left(U_1\cap \dots \cap U_n\right)$ is an open set.
\end{proof}

In view of Theorems \ref{Thm:Lewis-root} and \ref{Thm:main-root}, a root system is representable if and only if it is Esakia representable. Because of this, it is natural to ask whether every Priestley space whose underlying poset is a root system is also an Esakia space. The next example provides a negative answer to this question.

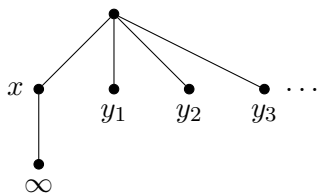
\begin{figure}
\[
	\begin{tikzpicture}
		\tikzstyle{point} = [shape=circle, thin, draw=black, fill=black, scale=0.35]
		\node (0) at (0,0)[point] {};
		\node[label=left:{$x$}] (x) at (-1,-1)[point] {};
		\node[label=below:{$\infty$}] (inf) at (-1,-2)[point] {};
		\node[label=below:{$y_1$}] (y1) at (0,-1)[point] {};
		\node[label=below:{$y_2$}] (y2) at (1,-1)[point] {};
		\node[label=below:{$y_3$}] (y3) at (2,-1)[point] {};
		\node (hinf) at (2.5,-1)[] {$\hdots$};
		\draw (0)--(x)--(inf);
		\draw (0)--(y1);
		\draw (0)--(y2);
		\draw (0)--(y3);
	\end{tikzpicture}
\]
\caption{An infinite root system which can be turned into a Priestley space that is not an Esakia space.}\label{Figure 2}
\end{figure}

\begin{exa}
Let $\langle X; \leq \rangle$ be the infinite root system depicted in Figure \ref{Figure 2}. When endowed with the topology
\[
\tau = \{ U \subseteq X : \text{either }\infty \notin U \text{ or }U\text{ is cofinite}\},
\]
the root system $\langle X; \leq \rangle$ becomes a Priestley space $\langle X; \leq, \tau\rangle$. We will show that $\langle X; \leq, \tau\rangle$ is not an Esakia space. Suppose the contrary, with a view to contradiction. Since $x$ is isolated, the downset ${\downarrow}x$ is open and, therefore, $X \smallsetminus {\downarrow}x$ is closed. As every point of $X$ other than $\infty$ is isolated, we obtain that  $X \smallsetminus {\downarrow}x$ is an infinite closed set whose members are all isolated points. Clearly, this contradicts the assumption that $\langle X; \leq, \tau\rangle$ is compact.
\qed
\end{exa}

\section{Esakia representable well-ordered forests}

Recall from Proposition \ref{Prop:closure}(\ref{Prop:closure:item:1}) that the class of representable posets is closed under order duals. Therefore, Theorem \ref{Thm:Lewis-root} can also be viewed as a characterization of the representable forests. More precisely, we have following.

\begin{Theorem}\label{Thm:Lewis-tree}
A forest is representable iff it has enough gaps and each of its nonempty chains has a supremum.
\end{Theorem}

It is therefore natural to wonder whether the above result holds for Esakia representable forests too. However, this is not the case because the tree depicted in Figure \ref{Figure 1} is not Esakia representable (see \cite[Example 5.6]{BeMo09}), although it has enough gaps and each of its nonempty chains has a supremum. Notice that the tree in Figure \ref{Figure 1} contains an \emph{infinite descending chain}
\[
\dots < x_n < \dots < x_2 < x_1 < x_0.
\]
Our main result states that the above description of the representable forests can be extended to Esakia representable forests by prohibiting the presence of such chains.

\begin{Definition}
A forest is \textit{well-ordered} when it lacks infinite descending chains, that is, it does not contain any subposet isomorphic to the order dual of $\langle \mathbb{N}; \leq \rangle$.\footnote{Well-ordered trees and forests are often endowed with a strict order relation. For the present purpose, however, it is convenient to endow them with a nonstrict order relation, so that they can be viewed as posets. The two presentations are of course equivalent.}
\end{Definition}

Notice that every well-ordered forest has enough gaps. Therefore, our main result takes the following form.

\begin{Theorem}\label{Thm:main-forest}
A well-ordered forest is Esakia representable iff each of its nonempty chains has a supremum.
\end{Theorem}

Let $X = \langle X; \leq \rangle$ be a well-ordered forest. We recall for each $x \in X$ there exists a unique ordinal $\alpha$ such that $\langle {\downarrow}x \smallsetminus \{ x \}; < \rangle$ is isomorphic to $\langle \alpha; \in \rangle$. The ordinal $\alpha$ is called the \emph{order type} of $x$ and will be denoted by $\h{x}$. Given $Y \sub X$ and an ordinal $\alpha$, we let
\begin{align*}
\mathsf{h}\left(X\right)
	&
		\coloneqq \text{ the least ordinal $\alpha$ s.t. $\h{x} \leq \alpha$ for every $x \in X$};\\
X_\alpha
	&
		\coloneqq \{ x \in X : \h{x} = \alpha \};\\
X_{\ast \alpha}
	&
		\coloneqq \{ x \in X : \h{x}\ast \alpha \} \text{ for } \ast \in \{\leq, <, \geq, >\};\\
\color{black}
{\uparrow_\alpha}Y
	&
		\coloneqq X_{\leq\alpha} \cap {\uparrow}Y.
\end{align*}

The implication from left to right in Theorem \ref{Thm:main-forest} holds by Proposition \ref{Prop:properties}. The rest of the paper is devoted to proving the implication from right to left. As in the case of Theorem \ref{Thm:main-root}, it suffices to prove this implication for well-ordered trees (as opposed to arbitrary well-ordered forests). Therefore, from now on we fix an arbitrary well-ordered tree $X = \langle X; \leq \rangle$ in which every nonempty chain has a supremum. Our aim is to prove that $X$ is Esakia representable. To this end, we will define a topology $\tau_\alpha$ on $X_{\leq \alpha}$ for each ordinal $\alpha$ and show that $\langle X; \leq, \tau_{\h{X}}\rangle$ is indeed an Esakia space (observe that $X = X_{\leq \h{X}}$). 

\begin{figure}
\[
    \begin{tikzpicture}[scale=0.4]
		\tikzstyle{point} = [shape=circle, thin, draw=black, fill=black, scale=0.35]

    % labelling the levels
        % X_{\alpha+1}
            \node at (-14, 2) {$X_{\alpha+1}$};
            
        % X_\alpha
            \node at (-14, 0) {$x, y \in X_\alpha$};

    % labelling the open sets

    % labelling U
                \node at (-9, 4.8)
                    {
                        $
                            \textcolor{blue}{U}
                                =  
                                    \left(\textcolor{ForestGreen}{V}
                                    \, \cup \, {\uparrow_{\alpha+1}}
                                    \left(\textcolor{ForestGreen}{V} \cap X_\alpha\right)\right)
                                    \smallsetminus
                                    {\downarrow} Z
                        $
                    };

    % DRAW THE OPEN SETS FIRST SO THAT THE POINTS ARE NOT OVERSHADOWED      

        % U = V \cup {\uparrow}\left(V \cap X_\alpha\right) = blue trapezoid
            
            % from x_0 till x_1
            
                \fill[blue!10, rounded corners=10pt, opacity=0.3]
                    (-9.5, -6) -- (-12, 4) -- (-1.7, 4) -- (-1, -6) -- cycle;
                \draw[blue, rounded corners=10pt]
                    (-9.5, -6) -- (-12, 4) -- (-1.7, 4) -- (-1, -6) -- cycle;

            % from x_3 onwards
                \fill[blue!10, rounded corners=10pt, opacity=0.3] 
                    (5.5, -6) rectangle (13, 1);
                \draw[blue, rounded corners=10pt] 
                    (5.5, -6) rectangle (13, 1);

        % V = green rectangle
        
            % labelling V
                \node at (6.5, 2) {$\textcolor{ForestGreen}{V} \in \tau_\alpha$};
            % from x_0 till x_1
                \fill[ForestGreen!20, rounded corners=8pt, opacity=0.5] 
                (-9, 0.6) rectangle (-1.8, -5);
                \draw[ForestGreen, rounded corners=8pt] 
                (-9, 0.6) rectangle (-1.8, -5);
            % from x_3 onwards
                \fill[ForestGreen!20, rounded corners=8pt, opacity=0.5] 
                (6, 0.5) rectangle (12, -5);
                \draw[ForestGreen, rounded corners=8pt] 
                (6, 0.5) rectangle (12, -5);

        % x_0
	   \node[label = right:{$x$}] (x0) at (-8,0)[point] {};
            % cutting out x0
                \draw[red, very thick] (x0) ++(-0.3, -0.3) -- ++(0.6, 0.6);
                \draw[red, very thick] (x0) ++(-0.3, 0.3) -- ++(0.6, -0.6);
        % successors of x_0
            \node (x00) at (-11,2)[shape=circle, thin, draw=BurntOrange, fill=BurntOrange, scale=0.6] {};
            \node (x01) at (-10,2)[point] {};  
            \node[label = above : {$z$}] (x02) at (-9,2)[point] {};
            % cutting out x02
                \draw[red, very thick] (x02) ++(-0.3, -0.3) -- ++(0.6, 0.6);
                \draw[red, very thick] (x02) ++(-0.3, 0.3) -- ++(0.6, -0.6);
            \node (1dots) at (-8,2)[] {$\dots$}; 
        % predecessors of x_0
            \node (x0-1) at (-8,-1)[point] {};
            \node (x0-2) at (-8,-2)[point] {};
            \node (x0-3) at (-8,-3)[point] {};
            \node (x0-4) at (-8,-4)[] {$\vdots$};
            % cutting out the downset of a nonspecial successor of x_0
                \draw[red, very thick] (x0-1) ++(-0.3, -0.3) -- ++(0.6, 0.6);
                \draw[red, very thick] (x0-1) ++(-0.3, 0.3) -- ++(0.6, -0.6);
                \draw[red, very thick] (x0-2) ++(-0.3, -0.3) -- ++(0.6, 0.6);
                \draw[red, very thick] (x0-2) ++(-0.3, 0.3) -- ++(0.6, -0.6);
                \draw[red, very thick] (x0-3) ++(-0.3, -0.3) -- ++(0.6, 0.6);
                \draw[red, very thick] (x0-3) ++(-0.3, 0.3) -- ++(0.6, -0.6);
                \draw[red, very thick] (x0-4) ++(-0.3, -0.3) -- ++(0.6, 0.6);
                \draw[red, very thick] (x0-4) ++(-0.3, 0.3) -- ++(0.6, -0.6);
        % connecting x_0 with its upsed and downset
            \draw (x0)--(x0-1)--(x0-2)--(x0-3);
            \draw (x0)--(x00);
            \draw (x0)--(x01); 
            \draw (x0)--(x02);
            
        % x_1
            \node[label=right:{$y$}] (x1) at (-3,0)[point] {};
        % successors of x_1
            \node (x10) at (-6,2)[shape=circle, thin, draw=BurntOrange, fill=BurntOrange, scale=0.6] {};
            \node (x11) at (-5,2)[point] {};
            \node (x12) at (-4,2)[point] {};
            \node (1dots) at (-3,2)[] {$\dots$};
        % predecessors of x_1
            \node (x1-1) at (-3,-1)[point] {};
            \node (x1-2) at (-3,-2)[point] {};
            \node (x1-3) at (-3,-3)[point] {};
            \node (x1-4) at (-3,-4)[] {$\vdots$};  
        % connecting x_0 with its upset and downset
		\draw (x1)--(x1-1)--(x1-2)--(x1-3);
            \draw (x1)--(x10);
            \draw (x1)--(x11); 
            \draw (x1)--(x12);
        % cutting out the downset of x_1
                \draw[red, very thick] (x1) ++(-0.3, -0.3) -- ++(0.6, 0.6);
                \draw[red, very thick] (x1) ++(-0.3, 0.3) -- ++(0.6, -0.6);
                \draw[red, very thick] (x1-1) ++(-0.3, -0.3) -- ++(0.6, 0.6);
                \draw[red, very thick] (x1-1) ++(-0.3, 0.3) -- ++(0.6, -0.6);
                \draw[red, very thick] (x1-2) ++(-0.3, -0.3) -- ++(0.6, 0.6);
                \draw[red, very thick] (x1-2) ++(-0.3, 0.3) -- ++(0.6, -0.6);
                \draw[red, very thick] (x1-3) ++(-0.3, -0.3) -- ++(0.6, 0.6);
                \draw[red, very thick] (x1-3) ++(-0.3, 0.3) -- ++(0.6, -0.6);
                \draw[red, very thick] (x1-4) ++(-0.3, -0.3) -- ++(0.6, 0.6);
                \draw[red, very thick] (x1-4) ++(-0.3, 0.3) -- ++(0.6, -0.6);

        % x_2
            \node (x2) at (2,0)[point] {};
        % successors of x_2
            \node (x20) at (-1,2)[shape=circle, thin, draw=BurntOrange, fill=BurntOrange, scale=0.6] {};
            \node (x21) at (0,2)[point] {};
            \node (x22) at (1,2)[point] {};
            \node (1dots) at (2,2)[] {$\dots$};
        % predecessors of x_2
            \node (x2-1) at (2,-1)[point] {};
            \node (x2-2) at (2,-2)[point] {};
            \node (x2-3) at (2,-3)[point] {};
            \node (x2-4) at (2,-4)[] {$\vdots$};
        % connecting x_2 with its upset and downset
  		\draw (x2)--(x2-1)--(x2-2)--(x2-3);
            \draw (x2)--(x20);
            \draw (x2)--(x21); 
            \draw (x2)--(x22);

        % x_3
            \node (x3) at (7,0)[point] {};
        % predecessors of x_3
            \node (x3-1) at (7,-1)[point] {};
            \node (x3-2) at (7,-2)[point] {};
            \node (x3-3) at (7,-3)[point] {};
            \node (x3-4) at (7,-4)[] {$\vdots$};
        % connecting x_3 with its downset
            \draw (x3)--(x3-1)--(x3-2)--(x3-3);

        % remaining points
            \node (hdots) at (11,0) {$\dots$};

        % labelling Z
            \node at (-13.5, -3) {$Z = \{y, z\}$};
    \end{tikzpicture}
\]
\caption{A member of $\mathcal{S}_{\alpha+1}$ of the form described in condition (\ref{item:subb:3}). For each $v \in X_\alpha$ we coloured in \textcolor{BurntOrange}{orange} the corresponding element $v^+$ of $X_{\alpha+1}$. 
Furthermore, we coloured in \textcolor{ForestGreen}{green} the set $\textcolor{ForestGreen}{V} \in \tau_\alpha$. Lastly, $Z = \{ y, z \}$ is a finite subset of $P_\alpha \cup S_{\alpha+1}$. Then the set $\textcolor{blue}{U}
                                =  
                                    \left(\textcolor{ForestGreen}{V}
                                    \, \cup \, {\uparrow_{\alpha+1}}
                                    \left(\textcolor{ForestGreen}{V} \cap X_\alpha\right)\right)
                                    \smallsetminus
                                    {\downarrow} Z
                        $
                         is obtained by considering the \textcolor{blue}{blue} shape and removing the elements crossed in \textcolor{red}{red} from it.}
\label{Figure 3}
\end{figure}
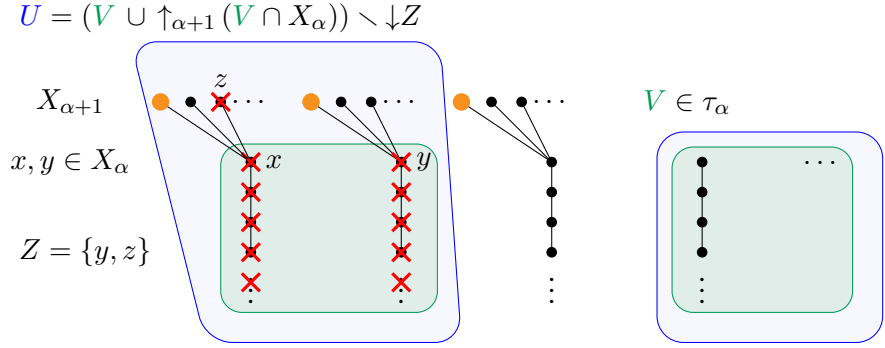

First, let $\tau_0$ be the unique topology on the singleton $X_{\leq 0}$. For the successor case, suppose that we already defined a topology $\tau_\alpha$ on $X_{\leq \alpha}$ for some ordinal $\alpha$. Then let
\[
P_\alpha \coloneqq \{ x \in X_{\alpha} : \exists y \in X_{\alpha +1} \text{ s.t. }x < y \}
\]
and for each $x \in P_\alpha$ choose an element $x^+ \in X_{\alpha +1}$ such that $x < x^+$.  Moreover, let
\[
S_{\alpha + 1}\coloneqq X_{\alpha + 1} \smallsetminus \{ x^+ : x \in P_\alpha \}.
\]
Lastly, let $\tau_{\alpha + 1}$ be the topology on $X_{\leq \alpha +1}$ generated by the subbase $\mathcal{S}_{\alpha+1}$ comprising the sets
\benroman
\item\label{item:subb:1} $\{ x \}$ for every $x \in S_{\alpha+1}$;
\item\label{item:subb:2} ${\downarrow} x$ for every $x \in P_\alpha$;
\item\label{item:subb:3} $\left(V \cup {\uparrow_{\alpha + 1}} \left(V \cap X_\alpha\right)\right) \smallsetminus {\downarrow} Z$ for every $V \in \tau_\alpha$ and every $Z \sub_\omega P_\alpha \cup S_{\alpha+1}$ (see Figure \ref{Figure 3}). 
\eroman

For the limit case, let $\alpha$ be a limit ordinal and suppose that we already defined a topology $\tau_\beta$ on $X_{\leq \beta}$ for each $\beta < \alpha$. Then let $\tau_\alpha$ be the topology on $X_{\leq \alpha}$ generated by the subbase
\[
\mathcal{S}_\alpha \coloneqq \{ V \cup {\uparrow_\alpha} \left(V \cap X_\beta\right) : \beta < \alpha \text{ and } V \in \tau_\beta \}.
\]

The next observation will be used later on.

\begin{Lemma}\label{lemma : preservation of opens}
For each pair of ordinals $\beta < \alpha$ and $U \in \tau_\beta$ we have $U \cup {\uparrow_\alpha}\left(U \cap X_\beta\right) \in \mathcal{S}_\alpha$.
\end{Lemma}

\begin{proof}
Let $U_\alpha \coloneqq U \cup {\uparrow_\alpha}\left(U \cap X_\beta\right)$. The proof proceeds by induction on $\alpha$. The case where $\alpha = 0$ holds vacuously because there exists no $\beta < 0$. For the successor case, we suppose that the statement holds for $\alpha$ and we will prove that it also holds for $\alpha +1$. Consider $\beta < \alpha +1$. We have two cases: either $\beta = \alpha$ or $\beta < \alpha$.

First, suppose that $\beta = \alpha$. Then $U \in \tau_\beta = \tau_\alpha$ by assumption. Therefore, condition (\ref{item:subb:3}) in the  definition of $\mathcal{S}_{\alpha + 1}$ and the assumption that $U \in \tau_\alpha$ guarantee that 
\[
U_{\alpha+1} = U \cup {\uparrow_{\alpha+1 }}\left(U \cap X_\alpha\right) \in \mathcal{S}_{\alpha +1}.
\]

Then we consider the case where $\beta < \alpha$.  By the inductive hypothesis we have $U_{\alpha} \in \tau_{\alpha}$. Thus, condition (\ref{item:subb:3}) in the  definition of $\mathcal{S}_{\alpha + 1}$ and guarantees that
\begin{equation}\label{Eq : first technical lemma section well ordered}
U_\alpha \cup {\uparrow_{\alpha +1 }}\left(U_\alpha \cap X_\alpha\right) \in \mathcal{S}_{\alpha + 1}.
\end{equation}
We claim that 
\begin{equation}\label{Eq : 2nd technical lemma section well ordered}
U_{\alpha+1} = U_\alpha \cup {\uparrow_{\alpha +1 }}\left(U_\alpha \cap X_\alpha\right).
\end{equation}
Together with condition (\ref{Eq : first technical lemma section well ordered}), this would imply $U_{\alpha+1} \in \mathcal{S}_{\alpha+1}$ as desired. 

To prove condition (\ref{Eq : 2nd technical lemma section well ordered}), consider $x \in U_{\alpha+1} = U \cup {\uparrow_{\alpha+1}}\left(U \cap X_\beta\right)$. If $x \in U$, then $x \in U_\alpha$ too by the definition of $U_\alpha$ and we are done.
Then we consider the case where $x \in {\uparrow_{\alpha+1}}\left(U \cap X_\beta\right)$. We have two cases: either $x \in {\uparrow_{\alpha}}\left(U \cap X_\beta\right)$ or $x \in X_{\alpha+1}$. If $x \in {\uparrow_{\alpha}}\left(U \cap X_\beta\right)$, then $x \in U_{\alpha}$ by the definition of $U_\alpha$ and we are done. Then we consider the case where $x \in X_{\alpha+1}$. Let $y$ be the unique member of $X_\alpha \cap {\downarrow}x$. Since $x \in {\uparrow_{\alpha+1}}\left(U \cap X_\beta\right)$ and $\beta \leq \alpha$, we have $y \in X_\alpha \cap {\uparrow_\alpha}\left(U \cap X_\beta\right)$. By the definition of $U_\alpha$ this yields $y \in U_\alpha \cap X_\alpha$ and, therefore, $x \in {\uparrow_{\alpha+1}}\left(U_\alpha\cap X_\alpha\right)$ as desired.

It only remains to prove the inclusion from right to left in condition (\ref{Eq : 2nd technical lemma section well ordered}). Consider $x \in U_\alpha \cup {\uparrow_{\alpha +1 }}\left(U_\alpha \cap X_\alpha\right)$. We have two cases: either $x \in U_\alpha$ or $x \in {\uparrow_{\alpha +1 }}\left(U_\alpha \cap X_\alpha\right)$. First, suppose that $x \in U_\alpha = U \cup {\uparrow_\alpha} \left(U \cap X_\beta\right)$. If $x\in U$, then $x \in U_{\alpha+1}$ by the definition of $U_{\alpha+1}$ and we are done. While if $x \in {\uparrow_\alpha} \left(U \cap X_\beta\right)$, then $x \in {\uparrow_{\alpha+1}}\left(U \cap X_\beta\right) \sub U_{\alpha+1}$ as desired. Then we consider the case where $x \in {\uparrow_{\alpha +1 }}\left(U_\alpha \cap X_\alpha\right)$. There exists $y \in U_\alpha \cap X_\alpha$ such that $y \leq x$. Since $U \in \tau_\beta$ by assumption, we have $U \subseteq X_{\leq \beta}$. Together with $\beta < \alpha$ and $y \in X_\alpha$, this yields $y \notin U$. Therefore, from $y \in U_\alpha = U \cup {\uparrow_\alpha}\left(U \cap X_\beta\right)$ it follows that there exists $z \in U \cap X_\beta$ such that $z \leq y$. As $y \leq x$ and $x \in X_{\leq \alpha+1}$, we conclude that $x \in {\uparrow_{\alpha+1}}\left(U \cap X_\beta\right)$. Hence, $x \in U_{\alpha+1}$ as desired. This establishes condition (\ref{Eq : 2nd technical lemma section well ordered}) and concludes the analysis of the successor case.

Lastly, consider the case where $\alpha$ is a limit ordinal. Since $\beta < \alpha$ and $U \in \tau_\beta$, the definition of $\mathcal{S}_\alpha$ ensures that $U_\alpha \in \mathcal{S}_\alpha$.
\end{proof}

We shall now define a function that will play an important role in the compactness proof. For every $x \in X$ and ordinal $\alpha \geq \h{x}$ we define and element $f_x\left(\alpha\right) \in X$ by recursion as
\begin{align*}
	f_x\left(\h{x}\right) & \coloneqq  x; \\
	\\
	f_x\left(\alpha +1\right) & \coloneqq  \begin{cases}
		\left(f_x\left(\alpha\right)\right)^{+} & \text{ if } f_x\left(\alpha\right) \in P_\alpha; \\
		f_x\left(\alpha\right) & \text{ otherwise};
	\end{cases} \\
	\\
	f_x\left(\alpha\right) & \coloneqq  \bigvee\{ f_x\left(\beta\right) : \h{x} \leq \beta < \alpha \}\text{ when $\alpha$ is a limit ordinal}.
\end{align*}

 Informally, we will regard $f_x$ as a function from $\{ \alpha \in \mathsf{Ord} : \h{x} \leq \alpha \}$ to $X$ (although its domain is not a set). Furthermore, given a pair of ordinals $\alpha$ and $\beta$, we write
\[
[\alpha, \beta]
	\coloneqq	
		\{ \gamma \in \mathsf{Ord} : \alpha \leq \gamma \leq \beta \}
		 \, \, \text{ and } \, \,
[\alpha, \beta) \coloneqq \{ \gamma \in \mathsf{Ord} : \alpha \leq \gamma < \beta \}.
\]

\begin{Lemma}\label{Lem : fx is well defined}
For every $x \in X$ the function $f_x$ is well defined and order preserving.
\end{Lemma}

\begin{proof}
It suffices to prove that for every ordinal $\alpha \geq \h{x}$ the restriction $f_x \colon [\h{x}, \alpha] \to X$ is well defined and order preserving. The proof works by induction starting at $\h{x}$. The base case and the successor case are straightforward. Then we consider the case where $\alpha$ is a limit ordinal such that $\h{x} < \alpha$. By the inductive hypothesis $f_x \colon [\h{x}, \alpha) \to X$ is well defined and order preserving. Consequently, $\{ f_x\left(\beta\right) : \h{x} \leq \beta < \alpha \}$ is a chain which, moreover, is nonempty because $\h{x} < \alpha$. Therefore, this chain has a supremum $f_x\left(\alpha\right)$ in $X$ by assumption. Hence, $f_x \colon [\h{x}, \alpha] \to X$ is also well defined and order preserving.
\end{proof}

We will make use of the following properties of the function $f_x$.

\begin{Lemma}\label{lemma : properties of f_x}
	The following conditions hold for every $x, y \in X$ and ordinal $\alpha \geq \h{x}$:
	\benroman
		\item\label{label 1x} $f_x\left(\alpha\right) \in \max X_{\leq \alpha}$;
		\item\label{label 2x} $f_x\left(\alpha +1\right) \notin S_{\alpha +1}$;
		\item\label{label 3x} for every $y \leq f_x\left(\alpha\right)$ such that $\h{x} \leq \h{y}$ we have $y = f_x\left(\h{y}\right)$;
		\item\label{label 4x} $\h{x} \leq \h{f_x\left(\alpha\right)}$ and $f_x\left(\alpha\right) = f_x\left(\h{f_x\left(\alpha\right)}\right)$;
		\item\label{label 5x} for every $\beta \in [\h{f_x\left(\alpha\right)},\alpha]$ we have $f_x\left(\alpha\right) = f_x\left(\beta\right)$.
	\eroman
\end{Lemma}
\begin{proof}
In this proof will make extensive use of the fact that $f_x$ is order preserving (see Lemma \ref{Lem : fx is well defined}).

A straightforward induction on $\alpha$ establishes condition (\ref{label 1x}).  Condition (\ref{label 2x}) follows from (\ref{label 1x}) and the definition of $f_x$.  To prove condition (\ref{label 3x}), assume that $y \leq f_x\left(\alpha\right)$ and $\h{x} \leq \h{y}$. We will prove that $\h{y} \leq \alpha$. Suppose, on the contrary, that $\alpha < \h{y}$. By condition (\ref{label 1x}) we have $f_x\left(\alpha\right) \in \max X_{\leq \alpha}$. Together with $\alpha < \h{y}$, this yields  $y \nleq f_x\left(\alpha\right)$, a contradiction. Since $\h{y} \leq \alpha$ and $\h{x} \leq \h{y}$, we obtain $f_x\left(\h{y}\right) \leq f_x\left(\alpha\right)$. On the other hand, $y \leq  f_x\left(\alpha\right)$ by assumption. Therefore, the elements $y$ and $f_x\left(\h{y}\right)$ are comparable because $X$ is a tree. By condition (\ref{label 1x}) we have $f_x\left(\h{y}\right) \in \max X_{\leq \h{y}}$. This yields  $y \nless f_x\left(\h{y}\right)$ and $f_x\left(\h{y}\right) \nless y$. As $y$ and $f_x\left(\h{y}\right)$ are comparable, we conclude that $y = f_x\left(\h{y}\right)$ as desired. Then we turn to prove condition (\ref{label 4x}).\ As $\h{x} \leq \alpha$, we also have $x = f_x\left(\h{x}\right) \leq f_x\left(\alpha\right)$, whence $\h{x} \leq \h{f_x\left(\alpha\right)}$. By applying condition (\ref{label 3x}) to $y \coloneqq f_x\left(\alpha\right)$ we obtain $f_x\left(\alpha\right) = f_x\left(\h{f_x\left(\alpha\right)}\right)$.  Lastly, condition (\ref{label 5x}) is an immediate consequence of condition (\ref{label 4x}) and the fact that $f_x$ is order preserving.\color{black}
\end{proof}

\begin{Corollary}\label{corollary : f x alpha+1 in U in B}
Let $x \in X$ and $\alpha$ an ordinal such that $\h{x} \leq \alpha+1$. If $f_x\left(\alpha+1\right) \in U$ for some $U \in \mathcal{S}_{\alpha+1}$, there exist $V \in \tau_\alpha$ and $Z \sub_\omega P_\alpha \cup S_{\alpha+1}$ such that
	\[
		U = \left(V \cup {\uparrow_{\alpha+1}}\left(V \cap X_{\alpha} \right)\right)\smallsetminus {\downarrow}Z .
	\]
\end{Corollary}
\begin{proof}
As $U$ is a member of $\mathcal{S}_{\alpha+1}$, it satisfies one of the conditions (\ref{item:subb:1})--(\ref{item:subb:3}) in the definition of $\mathcal{S}_{\alpha+1}$. If $U$ satisfies condition (\ref{item:subb:3}), we are done. Then suppose that $U$ does not satisfy condition (\ref{item:subb:3}), with a view to contradiction. In this case, $U$ satisfies either condition  (\ref{item:subb:1}) or condition (\ref{item:subb:2}). If $U$ satisfies condition (\ref{item:subb:1}), there exists $y \in S_{\alpha+1}$ such that $U = \{ y \}$. Hence, $f_x\left(\alpha+1\right) \in U = \{ y \}$ and, therefore, $f_x\left(\alpha+1\right) = y \in S_{\alpha + 1}$, a contradiction with Lemma \ref{lemma : properties of f_x}(\ref{label 2x}). On the other hand, if $U$ satisfies condition (\ref{item:subb:2}), there exists $y \in P_{\alpha}$ such that $U = {\downarrow}y$. Therefore, $f_x\left(\alpha+1\right) \in U = {\downarrow}y$. Since $y \in P_\alpha$, we have $y \notin \max X_{\leq \alpha+1}$, whence $f_x\left(\alpha+1\right) \notin X_{\leq \alpha+1}$, a contradiction with Lemma \ref{lemma : properties of f_x}(\ref{label 1x}).
\end{proof}

\section{The main lemma}

The next result plays a central role in the proof that the topological space $\langle X; \tau_{\h{X}}\rangle$ is compact.

\begin{Main Lemma}\label{lemma : fundamental lemma}
	Let  $x \in X$ and $\alpha$ be an ordinal such that $\h{x} \leq \alpha$. If $f_x\left(\alpha\right) \in U$ for some $U \in \mathcal{S}_\alpha$, there exist
\[
		v \leq x, \quad
		Y \sub_\omega X_{> \h{x}} \cap {\uparrow_\alpha} v, \quad \text{and} \quad
		Z \sub_{\omega} X_{< \alpha} \cap {\uparrow}v
\]
	such that
	$
		{\uparrow_\alpha} v \smallsetminus \left({\uparrow_\alpha} Y \cup {\downarrow}Z\right) \sub U
	$
	and $\h{v}$ is either zero or a successor ordinal.
\end{Main Lemma}
\begin{proof}
It holds that $\h{x} \leq \alpha$ by assumption.\
	We proceed by induction on the left subtraction $\alpha - \h{x}$, i.e., the only ordinal $\beta$ such that $\h{x} + \beta = \alpha$.
	
\subsection*{Base case}
In the base case, $\alpha - \h{x} = 0$ and, therefore, $\h{x} = \alpha$. Together with the definition of $f_x$, this yields $x = f_x\left(\h{x}\right)=f_x\left(\alpha\right)$. Consequently,  Lemma \ref{lemma : properties of f_x}(\ref{label 1x}) implies $x \in \max X_{\leq\alpha}$, whence ${\uparrow_\alpha} x = \{x\}$.
Suppose first that either $\h{x} = 0$ or $\h{x}$ is a successor ordinal.\ Letting  $v \coloneqq x$, $Y \coloneqq \emptyset$, and  $Z \coloneqq \emptyset$
 and using the assumption that $x = f_x\left(\alpha\right)  \in U$, we obtain
\[
{\uparrow_\alpha}v \smallsetminus \left({\uparrow_\alpha}Y  \cup {\downarrow}Z\right) 
=
{\uparrow_\alpha}x \smallsetminus  \emptyset
=
\{ x \} \smallsetminus \emptyset
=
\{ x \} \subseteq U
\]
and we are done.  Then we consider the case where $\alpha = \h{x}$ is a limit ordinal. As $U \in \mathcal{S}_\alpha$ by assumption, the definition of $\mathcal{S}_\alpha$ implies that  there exist $\beta < \alpha$ and $V \in \tau_\beta$ such that $U = V \cup {\uparrow_\alpha}\left(V \cap X_\beta\right)$. From $\beta < \alpha$ and $V \in \tau_\beta$ it follows that $V\cap X_\alpha = \emptyset$ (because $V \subseteq X_{\leq \beta}$). As $\h{x} = \alpha$, this yields $x \notin V$. Together with the assumptions that $x = f_x\left(\alpha\right) \in U$ and $U = V \cup {\uparrow_\alpha}\left(V \cap X_\beta\right)$, this implies $x  \in {\uparrow_\alpha}\left(V \cap X_\beta\right)$. Consequently, there exists $v^\ast \leq x$ such that $v^\ast \in V \cap X_\beta$. Therefore,
\[
{\uparrow_\alpha}v^\ast  \sub {\uparrow_\alpha}\left(V \cap X_\beta\right) \sub U.
\]
Now, recall that $\beta < \alpha$ and that $\alpha$ is a limit ordinal. Therefore, there exists a successor ordinal $\gamma$ such that $\beta \leq \gamma < \alpha$. Furthermore, as $\h{x} = \alpha$, $\h{v^\ast} = \beta$, and $v^\ast \leq x$, there exists $v \in X$ such that $v^\ast \leq v \leq x$ and $\h{v} = \gamma$. In view of the above display and $v^\ast\leq v$, by letting $Y \coloneqq \emptyset$ and $Z \coloneqq \emptyset$, we conclude that
\[
{\uparrow_\alpha}v \smallsetminus \left({\uparrow_\alpha}Y \color{black} \cup {\downarrow}Z\color{black} \right) = {\uparrow_\alpha}v \smallsetminus \emptyset = {\uparrow_\alpha}v \subseteq  
   {\uparrow_\alpha} v^\ast \subseteq U.
\]
As $\h{v} = \gamma$ is a successor ordinal, we are done.

\subsection*{Successor case}
In the successor case of the induction, $\alpha - \h{x}$ is a successor ordinal $\beta +1$ and $\alpha = \h{x} + \beta +1$. By assumption we have
\[
f_x\left(\h{x} + \beta +1\right) = f_x\left(\alpha\right) \in U \, \, \text{ and } \, \, U \in \mathcal{S}_{\alpha} = \mathcal{S}_{\h{x} + \beta +1}.
\]
Therefore, we can apply Corollary \ref{corollary : f x alpha+1 in U in B} obtaining
	\begin{equation}\label{Eq : what U in the successor case is}
		U = \left(V \cup {\uparrow_{\alpha}}\left(V \cap X_{\h{x} + \beta} \right)\right)\smallsetminus {\downarrow}\bar{Z}
	\end{equation}
	for some $V \in \tau_{\h{x} + \beta}$ and $\bar{Z} \sub_\omega P_{\h{x} + \beta} \cup S_{\alpha}$.

\begin{Claim}\label{Claim : well-orderes claim 1}
$f_x\left(\h{x}+\beta\right) \in V$.
\end{Claim}

\begin{proof}[Proof of the Claim.]
Recall that $f_x\left(\alpha\right) \in U \subseteq V \cup {\uparrow_{\alpha}}\left(V \cap X_{\h{x} + \beta} \right)$. Therefore, we have two cases: either $f_x\left(\alpha\right) \in V$ or $f_x\left(\alpha\right)\in {\uparrow_{\alpha}}\left(V \cap X_{\h{x} + \beta} \right)$.\ First, suppose that $f_x\left(\alpha\right) \in V$. Then
\[
f_x\left(\alpha\right) \in V \subseteq X_{\leq \h{x} + \beta},
\]
where the last inclusion holds because $V \in \tau_{\h{x} + \beta}$. From from the above display and $\h{x} + \beta < \alpha$ it follows that $\h{f_x\left(\alpha\right)} \leq \h{x}+ \beta< \alpha$. By Lemma \ref{lemma : properties of f_x}(\ref{label 5x}) we conclude that $f_x\left(\h{x} + \beta\right) = f_x\left(\alpha\right) \in V$ as desired. Then we consider the case where $f_x\left(\alpha\right) \in {\uparrow_{\alpha}}\left(V \cap X_{\h{x} + \beta} \right)$. There exists $y \in V\cap X_{\h{x}+\beta}$ such that $y \leq f_x\left(\alpha\right)$. Since $y \leq f_x\left(\alpha\right)$ and $\h{x} \leq \h{y}$, we can apply  
	Lemma \ref{lemma : properties of f_x}(\ref{label 3x}), obtaining $y = f_x\left(\h{y}\right)$. As $y \in X_{\h{x}+\beta}$ and, therefore, $\h{y} = \h{x}+\beta$, we conclude that $f_x\left(\h{x}+\beta\right) = f_x\left(\h{y}\right) = y  \in V$.
\end{proof}

Now, recall that $\mathcal{S}_{\h{x}+\beta}$ is a subbase for $\tau_{\h{x}+\beta}$ and that $V \in \tau_{\h{x}+\beta}$. As $f_x\left(\h{x}+\beta\right) 	\in V$ by Claim \ref{Claim : well-orderes claim 1}, there exist $W_1, \dots, W_n \in \mathcal{S}_{\h{x}+\beta}$ such that
\begin{equation}\label{Eq : main lemma : successor : 2nd claim : IH}
	f_x\left(\h{x}+\beta\right) \in W_1 \cap \dots \cap W_n \subseteq V.
\end{equation}

\begin{Claim}\label{Claim : main lemma : 2nd}
There exist $v \leq x$, 		$Y^* \sub_\omega X_{> \h{x}} \cap {\uparrow_{\h{x}+\beta}} v$, and $Z^* \sub_\omega X_{< \h{x}+\beta} \cap {\uparrow} v$ such that 
\[
	{\uparrow_{\h{x}+\beta}} v
		\smallsetminus
			\left(
				{\uparrow_{\h{x}+\beta}} Y^*
					\cup
				{\downarrow} Z^*
			\right)
				\sub
					W_1 \cap \dots \cap W_n
						\subseteq
							V
\]
and $\h{v}$ is either zero or a successor ordinal.
\end{Claim}

\begin{proof}[Proof of the Claim.]
By applying the inductive hypothesis to $W_1, \dots, W_n \in \mathcal{S}_{\h{x}+\beta}$ and condition (\ref{Eq : main lemma : successor : 2nd claim : IH}), we obtain that for every $m \leq n$ \color{black} there exist
\[
	v_m \leq x,
		\quad
	Y_m \sub_\omega X_{> \h{x}} \cap {\uparrow_{\h{x}+\beta}} v_m, \, \, \text{ and } \, \, 
	Z_m \sub_\omega X_{<\h{x}+\beta} \cap {\uparrow} v_m
\]
such that 
\begin{equation}\label{Eq : ym has the desired property}
	{\uparrow_{\h{x}+\beta}} v_m
		\smallsetminus
			\left( {\uparrow_{\h{x}+\beta}} Y_m \cup {\downarrow} Z_m\right)\sub W_m
\end{equation}
and $\h{y_m}$ is either zero or a successor ordinal.\	As $X$ is a tree and $v_1, \dots, v_n \leq x$, the set $\{ v_m : m \leq n \}$ is a nonempty chain and, therefore, has a maximum $v$. Then, \color{black} letting
	\[
		Y^* \coloneqq \left(Y_1 \cup \dots \cup Y_n\right) \cap {\uparrow}v \, \, \text{ and } \, \, 
		Z^* \coloneqq \left(Z_1 \cup \dots \cup Z_m\right) \cap {\uparrow}v,
	\]
	we obtain
	\[
		Y^* \sub_\omega X_{> \h{x}} \cap {\uparrow_{\h{x}+\beta}} v
		\, \, \text{ and } \, \, 
		Z^* \sub_\omega X_{< \h{x}+\beta} \cap {\uparrow} v.
	\]
Furthermore, $v \leq x$ and $\h{v}$ is either zero or a successor ordinal. Therefore, it only remains to prove that
\[
	{\uparrow_{\h{x}+\beta}} v
		\smallsetminus
			\left(
				{\uparrow_{\h{x}+\beta}} Y^*
					\cup
				{\downarrow} Z^*
			\right)
				\sub
					W_1 \cap \dots \cap W_n
						\subseteq
							V.
\]

Since $W_1 \cap \dots \cap W_n
						\subseteq
							V$ by condition (\ref{Eq : main lemma : successor : 2nd claim : IH}), it suffices to show that ${\uparrow_{\h{x}+\beta}} v
		\smallsetminus
			\left(
				{\uparrow_{\h{x}+\beta}} Y^*
					\cup
				{\downarrow} Z^*
			\right)
				\sub
					W_1 \cap \dots \cap W_n$. To this end, consider $z \in {\uparrow_{\h{x}+\beta}} v
		\smallsetminus
		\left(
				{\uparrow_{\h{x}+\beta}} Y^*
					\cup
				{\downarrow} Z^*
		\right)$ and $m \leq n$. We need to show that $z \in W_m$. In view of condition (\ref{Eq : ym has the desired property}), the definition of $Y^\ast$ and $Z^\ast$, and $v_m \leq v$, it will be enough to show that $z \notin {\uparrow_{\h{x}+\beta}} Y_m\cup{\downarrow} Z_m$. We begin by proving that $z \notin {\uparrow_{\h{x}+\beta}} Y_m$. Suppose the contrary, with a view to contradiction. Then $\h{z} \leq \h{x} + \beta$ and there exists $y \in Y_m$ such that $y \leq z$. Since $v \leq z$ and $X$ is a tree, the elements $v$ and $y$ must be comparable. We have two cases: either $v \leq y$ or $y < v$. If $v \leq y$, then $y \in Y^*$ because $y \in Y_m$. Therefore, $z \in {\uparrow_{\h{x}+\beta}} Y^*$, a contradiction with $z \in {\uparrow_{\h{x}+\beta}} v
		\smallsetminus
			\left(
				{\uparrow_{\h{x}+\beta}} Y^*
					\cup
				{\downarrow} Z^*
			\right)$.
Then we consider the case where $y < v$. As $v \leq x$, this implies $\h{y} < \h{x}$, a contradiction with $y \in Y_m \subseteq X_{> \h{x}}$. Hence, we conclude that $z \notin {\uparrow_{\h{x}+\beta}} Y_m$. Then we turn to prove that $z \notin {\downarrow} Z_m$. Suppose the contrary, with a view to contradiction. Then there exists $y \in Z_m$ such that $z \leq y$. Since $v \leq z$, we obtain $v  \leq y \in Z_m$. By the definition of $Z^*$ we obtain $y \in Z^*$ and, therefore, $z \in {\downarrow}Z^*$, a contradiction with $z \in {\uparrow_{\h{x}+\beta}} v
		\smallsetminus
			\left(
				{\uparrow_{\h{x}+\beta}} Y^*
					\cup
				{\downarrow} Z^*
			\right)$. This establishes the above display.
\end{proof} 			
		
Now, consider the sets
\color{black}
	\[
		Y	\coloneqq
				Y^*
					\cup
				\left(
						  X_{\alpha} \cap \bar{Z}
						\cap 
					{\uparrow}v
				\right) \, \, \text{ and } \, \, 
		Z
			\coloneqq 
				Z^*
					\cup
				\left(
					 X_{\h{x}+\beta} \cap {\downarrow}\bar{Z}
						\cap
					{\uparrow}v
				\right).
	\]
From Claim \ref{Claim : main lemma : 2nd} it follows that $v \leq x$ and that $\h{v}$ is either zero or a successor ordinal. Furthermore, as $\bar{Z}$ and $Y^\ast$ are finite (the latter by Claim \ref{Claim : main lemma : 2nd}), the set $Y$ is also finite. Lastly, as $X$ is a tree and $\bar{Z}$ is finite, ${\downarrow}\bar{Z}$ is a union of finitely many chains. Therefore, $X_{\h{x}+\beta} \cap {\downarrow}\bar{Z}$ is a finite set. As $Z^\ast$ is finite by Claim \ref{Claim : main lemma : 2nd}, we conclude that $Z$ is also finite. Therefore, it only remains to show that 
\[
	Y \sub X_{> \h{x}} \cap {\uparrow_\alpha} v,
		\quad
	Z \sub X_{< \alpha} \cap {\uparrow}v,
	\, \, \text{ and } \, \, 
	{\uparrow_\alpha} v \smallsetminus
		\left(
			{\uparrow_\alpha} Y
				\cup
			{\downarrow} Z
		\right)
			\sub U.
\]

By Claim \ref{Claim : main lemma : 2nd} we have $Y^* \subseteq X_{> \h{x}} \cap {\uparrow_{\h{x}+\beta}} v \subseteq X_{> \h{x}} \cap {\uparrow_{\alpha}} v$ and $Z^* \subseteq X_{< \h{x}+\beta} \cap {\uparrow} v \subseteq X_{< \alpha} \cap {\uparrow} v$. Together with $\alpha = \h{x} + \beta + 1$ and the definition of $Y$ and $Z$, this guarantees the validity of the first two conditions in the above display. Therefore, it only remains to prove that ${\uparrow_\alpha} v \smallsetminus
		\left(
			{\uparrow_\alpha} Y
				\cup
			{\downarrow} Z
		\right)
			\sub U$. By condition (\ref{Eq : what U in the successor case is}) this amounts to
\begin{equation}\label{Eq : successor case - what we need to prove (final)}
	{\uparrow_\alpha} v \smallsetminus
		\left(
			{\uparrow_\alpha} Y
				\cup
			{\downarrow} Z
		\right)
			\sub
			\left(V \cup {\uparrow_{\alpha}}\left(V \cap X_{\h{x} + \beta} \right)\right)\smallsetminus {\downarrow}\bar{Z}.		
\end{equation}

\color{black} Consider $z \in {\uparrow_\alpha} v \smallsetminus
		\left(
			{\uparrow_\alpha} Y
				\cup
			{\downarrow} Z
		\right)$. Then $v \leq z \in X_{\leq \alpha}$ and $z \notin {\uparrow_\alpha} Y \cup {\downarrow}Z$. Since $\alpha = \h{x} + \beta +1$ and $z \in X_{\leq \alpha}$,  we have two cases: either $z \in X_{\leq \h{x}+\beta}$ or $z \in X_\alpha$. First, suppose that $z \in X_{\leq \h{x}+\beta}$.
Then
\[
	z \in
		\left(X_{\leq \h{x}+\beta}
			\cap
			{\uparrow_\alpha}v\right)
				\smallsetminus
			\left(
					{\uparrow_\alpha}Y \cup {\downarrow}Z
			\right)	
	\sub
		{\uparrow_{\h{x}+\beta}}v
			\smallsetminus
		\left(
			{\uparrow_{\h{x}+\beta}} Y^* \cup {\downarrow} Z^*
		\right)
	\sub V,
\]
where the first inclusion holds because $Y^* \subseteq Y$ and $Z^* \sub Z$, and the last by Claim \ref{Claim : main lemma : 2nd}. Therefore, in order to conclude that $z$ belongs to the right hand side of condition (\ref{Eq : successor case - what we need to prove (final)}), it suffices to show that $z \notin {\downarrow}\bar{Z}$.  Suppose the contrary, with a view to contradiction. Then there exists $y \in \bar{Z}$ such that $z \leq y$. Since $v \leq z \in X_{\leq \h{x} + \beta}$ and $\bar{Z} \subseteq P_{\h{x} + \beta} \cup S_{\alpha}$, there exists $y^\ast \in X_{\h{x} + \beta}$ such that $v \leq z \leq y^\ast \leq y$. Hence, $y^\ast \in  X_{\h{x} + \beta} \cap {\downarrow} \bar{Z} \cap {\uparrow}v \subseteq Z$, where the last inclusion holds by the definition of $Z$. Together with $z \leq y^\ast$, this implies $z \in {\downarrow}Z$, which is false.

Then we consider the case where $z \in X_\alpha$. Let $y$ be the unique element of $X_{\h{x}+\beta} \cap {\downarrow}z$. We will prove that $y \in V$.  By Claim \ref{Claim : main lemma : 2nd} it suffices to show that
\[
y \in {\uparrow_{\h{x}+\beta}} v
		\smallsetminus
			\left(
				{\uparrow_{\h{x}+\beta}} Y^*
					\cup
				{\downarrow} Z^*
			\right).
\]
 From $z \in X_\alpha$, $v \leq x$, and $\h{x} < \alpha$ it follows that $\h{v} < \h{z}$. Since $v \leq z$, this implies $v < z$. Moreover, as $X$ is a tree, from $v, y < z$ it follows that $v$ and $y$ are comparable. Since $y$ is the unique immediate predecessor of $z$ by definition and $v < z$, we conclude that $v \leq y$.  Hence, $y \in {\uparrow_{\h{x}+\beta}} v$. Now, observe that $y \notin {\uparrow} Y^*$,  otherwise we would have $z \in {\uparrow_\alpha}Y$, a contradiction. Moreover, observe that $y \notin {\downarrow}Z^*$ because $y \in X_{\h{x}+\beta}$ and $Z^* \sub X_{< \h{x}+\beta}$   by Claim \ref{Claim : main lemma : 2nd}. This establishes the above display and, therefore, that $y \in V$. Together with $y \leq z$, $y \in X_{\h{x}+ \beta}$, and $z \in X_\alpha$,  this yields $z \in {\uparrow_\alpha}\left(V \cap X_{\h{x}+\beta}\right)$.  Therefore, in order to prove that $z$ belongs to the right hand side of condition (\ref{Eq : successor case - what we need to prove (final)}),   it only remains to show that $z \notin {\downarrow}\bar{Z}$. Suppose the contrary, with a view of contradiction. Then there exists $u \in \bar{Z}$ such that $z \leq u$. Together with $z \in X_\alpha$ and $u \in \bar{Z} \sub P_{\h{x}+\beta} \cup S_\alpha$,  this implies $z = u \in \bar{Z}$. Hence, $z \in X_\alpha \cap \bar{Z} \cap {\uparrow}v$. By the definition of $Y$ this yields $z \in {\uparrow_\alpha}Y$, a contradiction.

\subsection*{Limit case} Finally, we consider the case where $\alpha -\h{x}$ is a limit ordinal. In this case, $\alpha$ is also a limit ordinal. Consequently, from $U \in \mathcal{S}_\alpha$ it follows that $U = V  \cup  {\uparrow_\alpha} \left(V \cap X_\beta\right)$ for some $\beta < \alpha$ and $V \in \tau_\beta$. We have two cases: either $\beta < \h{x}$ or $\h{x} \leq \beta$.
	
First, suppose that $\beta < \h{x}$. Then $\beta < \h{x} \leq \h{f_x\left(\alpha\right)}$. Together with $V \subseteq X_{\leq \beta}$ (because $V \in \tau_\beta$), this implies $f_x\left(\alpha\right) \notin V$. On the other hand, $f_x\left(\alpha\right) \in U = V \cup {\uparrow_\alpha} \left(V \cap X_\beta\right)$
by assumption. Therefore, $f_x\left(\alpha\right) \in {\uparrow_\alpha} \left(V \cap X_\beta\right)$. Then there exists  $z \in V \cap X_\beta$ such that $z \leq f_x\left(\alpha\right)$. As $X$ is a tree, from $x, z \leq f_x\left(\alpha\right)$ it follows that $x$ and $z$ are comparable. Since $z \in X_\beta$ and $\beta < \h{x}$, we deduce that $z < x$. Thus,
\[
	{\uparrow_\alpha}x
		\sub
	{\uparrow_\alpha} z
		\sub
	{\uparrow_\alpha}\left(V \cap X_\beta\right)
		\sub
	U. 
\]
Now, let $Y \coloneqq \emptyset$ and $Z \coloneqq \emptyset$. Furthermore, if $\h{x}$ is  zero or a successor ordinal, let $v \coloneqq x$. While if $\h{x}$ is a limit ordinal, recall that $z < x$ and let $v$ be any element strictly between $z$ and $x$ whose height is a successor ordinal. In both cases, we are done.
	
Then we consider the case where $\h{x} \leq \beta$.\ We will prove that $f_x\left(\beta\right) \in V$.\
Recall that $f_x\left(\alpha\right) \in U = V \cup {\uparrow_\alpha}\left(V \cap X_\beta\right)$.\ Then we have two cases: either $f_x\left(\alpha\right) \in V$ or $f_x\left(\alpha\right) \in {\uparrow_\alpha}\left(V \cap X_\beta\right)$. If $f_x\left(\alpha\right) \in V$, from $V \subseteq X_{ \leq \beta}$ it follows that $\h{f_x\left(\alpha\right)} \leq \beta$. Together with $\beta \leq \alpha$ and Lemma \ref{lemma : properties of f_x}(\ref{label 5x}), this yields $f_x\left(\beta\right) = f_x\left(\alpha\right) \in V$ as desired. Then we consider the case where $f_x\left(\alpha\right) \in {\uparrow_\alpha}\left(V \cap X_\beta\right)$.  There exists $z \in V$ such that $\h{z} = \beta$ and $z \leq f_x\left(\alpha\right)$.  By Lemma \ref{Lem : fx is well defined}(\ref{label 3x}) we have $f_x\left(\beta\right) = z \in V$.\ This establishes that $f_x\left(\beta\right) \in V$ as desired.
	
As $\mathcal{S}_\beta$ is a subbase for the topology $\tau_\beta$, from $f_x\left(\beta\right) \in V \in \tau_\beta$ it follows that there exist $W_1, \dots, W_n \in \mathcal{S}_\beta$ such that $f_x\left(\beta\right) \in W_1 \cap \dots \cap W_n$. 
 Since $\h{x} \leq \beta < \alpha$, we have $\beta - \h{x} < \alpha - \h{x}$. Therefore, we can apply the inductive hypothesis obtaining  that for each $m \leq n$ there exist
\[
	v_m \leq x,
		\quad
	Y_m \sub_\omega X_{> \h{x}} \cap {\uparrow_{\beta}} v_m,
		\quad
	Z_m \sub_\omega X_{<\beta} \cap {\uparrow} v_m
\]
such that
\[
	{\uparrow_{\beta}} v_m
		\smallsetminus
			\left( {\uparrow_{\beta}} Y_m \cup {\downarrow} Z_m\right)\sub W_m
\]
and $\h{v_m}$ is either zero or a successor ordinal.\ As $X$ is a tree and $v_1, \dots, v_n \leq x$, the set $\{ v_m : m \leq n \}$ is a nonempty chain and, therefore, has a maximum $v$.
Then, letting
	\[
		Y \coloneqq \left(Y_1 \cup \dots \cup Y_n\right) \cap {\uparrow}v \, \, \text{ and } \, \, 
		Z \coloneqq \left(Z_1 \cup \dots \cup Z_m\right) \cap {\uparrow}v,
	\]
	we obtain
\begin{equation}\label{Eq : Y tilde-a}
		Y \sub_\omega X_{> \h{x}} \cap {\uparrow_{\beta}} v\, \, \text{ and } \, \,
		Z \sub_\omega X_{<\beta} \cap {\uparrow} v.
\end{equation}
Furthermore,
\begin{equation}\label{Eq : Y tilde}
	{\uparrow_{\beta}} v
		\smallsetminus
			\left(
				{\uparrow_{\beta}} Y
					\cup
				{\downarrow} Z
			\right)
				\sub
					W_1 \cap \dots \cap W_n
						\subseteq
							V
\end{equation}
and $\h{v}$ is either zero or a successor ordinal.

From condition (\ref{Eq : Y tilde-a}) and $\beta \leq \alpha$ it follows that
\[
Y \sub_\omega X_{> \h{x}} \cap {\uparrow_{\alpha}} v\, \, \text{ and } \, \,
		Z \sub_\omega X_{<\alpha} \cap {\uparrow} v.
\]
Since $\h{v}$ is either zero or a successor ordinal and $U = V \cup {\uparrow_\alpha}\left(V \cap X_\beta\right)$, it only remains to show that 
\[
	{\uparrow_\alpha} v \smallsetminus \left({\uparrow_\alpha} Y \cup {\downarrow}Z\right) \sub V \cup {\uparrow_\alpha}\left(V \cap X_\beta\right). 
\]
To this end, let $z \in X_{\leq \alpha}$ be such that $v \leq z$ and $z \notin {\uparrow_\alpha}Y \cup {\downarrow}Z$. We have two cases: either $z \in X_{\leq \beta}$ or $z \notin X_{\leq \beta}$. In the former case, we have $z \in {\uparrow_{\beta}} v
		\smallsetminus
			\left(
				{\uparrow_{\beta}} Y
					\cup
				{\downarrow} Z
			\right)$. By condition (\ref{Eq : Y tilde}) we conclude that $z \in V$ as desired. Then we consider the case where $z \notin X_{\leq \beta}$, i.e., $\h{z} > \beta$. Let $y$ be the unique element of ${\downarrow}z \cap X_\beta$. As  $X$ is a tree and $v, y \leq z$, we deduce that either $y < v$ or $v \leq y$. However, the former case cannot happen because $y \in X_\beta$, $v \leq x$, and $\h{x} \leq \beta$. Hence, $v \leq y$ and, therefore, $y \in {\uparrow_\beta}v$. Moreover, $y \notin {\uparrow_\beta}Y$ because $z \notin {\uparrow_\alpha}Y$ and $y \leq z \in X_{\leq \alpha}$. Lastly, $y \notin {\downarrow}Z$ because $y \in X_\beta$ and $Z \sub X_{< \beta}$. Therefore, $y \in {\uparrow_{\beta}} v
		\smallsetminus
			\left(
				{\uparrow_{\beta}} Y
					\cup
				{\downarrow} Z
			\right)$. From condition (\ref{Eq : Y tilde}) it follows that $y \in V$. Thus, $y \in V \cap X_\beta$. Together with $y \leq z \in X_{\leq \alpha}$, this implies $z \in {\uparrow_\alpha}\left(V \cap X_\beta\right)$.
\end{proof}

\section{Compactness}

The aim of this section is to prove the following.

\begin{Theorem}\label{Thm : Esakia compactness : trees}
The topological space $\langle X; \tau_{\h{X}} \rangle$ is compact.
\end{Theorem}

\noindent To this end, let $\mathcal{C}$ be an open covering of $X$. We need to show that $\mathcal{C}$ has a finite subcover. By Alexander's subbase theorem we may assume $\mathcal{C} \sub \mathcal{S}_{\h{X}}$.\ The construction of a finite subcover proceeds through a series of technical observations.

\begin{Proposition}\label{proposition : covering X minus down Z suffices}
	For every $x \in X$ there exists $\mathcal{V}_x \sub_\omega \mathcal{C}$ such that ${\downarrow}x \sub \bigcup \mathcal{V}_x$. 
\end{Proposition}
\begin{proof}
	We proceed by induction on $\h{x}$.
	If $\h{x} = 0$, then $x$ is the root of $X$ and the claim follows from ${\downarrow}x = \{x \}$. 
 If $\h{x} = \alpha+1$, there exists $y \in X_\alpha$ such that ${\downarrow}x = \{x\} \cup {\downarrow}y$. As $\h{y} = \alpha < \alpha +1$, by the inductive hypothesis there exists $\mathcal{V}_y \sub_\omega \mathcal{C}$ such that ${\downarrow}y \sub \bigcup \mathcal{V}_y$. Let $U \in \mathcal{C}$ be such that $x \in U$. Letting $\mathcal{V}_x \coloneqq \mathcal{V}_y \cup \{ U \}$, we conclude that $\mathcal{V}_x \sub_\omega \mathcal{C}$ and ${\downarrow}x \sub \bigcup \mathcal{V}_x$.\color{black}

Finally, suppose that $\h{x}$ is a limit ordinal and consider $U \in \mathcal{C}$ such that $x \in U$. We begin with the following observation.

\begin{Claim}\label{claim : y < x}
There exists $y < x$ such that $[y,x] \sub U$. 
\end{Claim}

\begin{proof}[Proof of the Claim]
We will prove that for every $\alpha \leq \h{X}$,
\[
\text{if }x \in W \text{ for some $W \in \mathcal{S}_\alpha$, there exists }y < x\text{ such that }[y,x] \sub W. 
\]
Since $x \in U \in \mathcal{C} \subseteq \mathcal{S}_{\h{X}}$ by assumption, Claim \ref{claim : y < x} follows immediately from the above display in the case where $\alpha = \h{X}$ and $W = U$.

We proceed by induction on $\alpha$. The case where $\alpha = 0$ is straightforward because the assumption that $\h{x}$ is a limit ordinal guarantees that $x \notin X_{\leq 0}$ and, therefore, $x \notin \bigcup \mathcal{S}_{0}$. Then we consider the case where $\alpha$ is a successor ordinal $\beta + 1$. Suppose that $x \in W \in \mathcal{S}_{\beta +1}$. Since $\h{x}$ is a limit ordinal, we have $x \notin X_{\beta +1}$. Therefore, the definition of $\mathcal{S}_{\beta+1}$ and $x \in W \in \mathcal{S}_{\beta+1}$ ensures that either $W = {\downarrow}z$ for some $z \in P_\beta$ or $W = V \cup {\uparrow_{\leq \beta+1}}\left(V \cap X_\beta\right)$ for some $V \in \tau_\beta$. First, suppose that $W = {\downarrow}z$. As $\h{x}$ is a limit ordinal, there exists $y < x$. Since $x \in W = {\downarrow}z$, we obtain $[y, x] \subseteq W$ as desired. Then we consider the case where $W = V \cup {\uparrow_{\leq \beta+1}}\left(V \cap X_\beta\right)$ for some $V \in \tau_\beta$. Together with $x \in W \smallsetminus  X_{\beta +1}$, this yields $x \in V$. As $V \in \tau_\beta$ and $\mathcal{S}_\beta$ is a subbase for $\tau_\beta$, there exist $V_1,\dots,V_n \in  \mathcal{S}_\beta$ such that $x \in V_1 \cap \dots \cap V_n = V \sub W$. Then the inductive hypothesis ensures that there exist $y_1, \dots, y_n < x$ such that $[y_1,x] \sub V_1, \dots, [y_n,x] \sub V_n$. As $y_1, \dots, y_n < x$ and $X$ is a tree, the set $\{ y_1, \dots, y_n \}$ is a nonempty chain and, therefore, has a maximum $y$. We have $y < x$ and $[y,x] \sub V_1 \cap \dots \cap V_n \sub W$ as desired.

Lastly, we consider the case where $\alpha$ is a limit ordinal.\ Suppose that $x \in W \in \mathcal{S}_{\alpha}$. Then $W = V \cup {\uparrow}\left(V \cap X_\beta\right)$ for some $\beta < \alpha$ and $V \in \tau_\beta$ by the definition of $\mathcal{S}_{\alpha}$. Together with $x \in W$, this yields $x \in V \cup {\uparrow}\left(V \cap X_\beta\right)$. If $x \in {\uparrow}\left(V \cap X_\beta\right)$, there exists $y \in V \cap X_\beta$ such that $y < x$. Therefore, $[y,x] \sub {\uparrow}\left(V \cap X_\beta\right) \sub W$ and we are done. Then we consider the case where $x \in V$. Together with $V \subseteq W$, the assumption that $V \in \tau_\beta$ and that $\mathcal{S}_\beta$ is a subbase for $\tau_\beta$ implies the existence of $V_1,\dots,V_n \in  \mathcal{S}_\beta$ such that $x \in V_1 \cap \dots \cap V_n = V \sub W$. Since $\beta < \alpha$, we can apply the inductive hypothesis obtaining $y_1, \dots, y_n < x$ such that $[y_1,x] \sub V_1, \dots, [y_n,x] \sub V_n$. As before, letting $y$ be the maximum of $\{ y_1, \dots, y_n \}$, we obtain $y < x$ and $[y,x] \sub V_1 \cap \dots \cap V_n \sub W$. 
\end{proof}

By the Claim there exists $y < x$ such that $[y, x] \subseteq U$.\ As $\h{y} < \h{x}$, we can apply the inductive hypothesis, obtaining that there exists $\mathcal{V}_y \sub_\omega \mathcal{C}$ such that ${\downarrow}y \sub \bigcup \mathcal{V}_y$. Therefore, letting $\mathcal{V}_x \coloneqq \mathcal{V}_y \cup \{ U \}$, we obtain that $\mathcal{V}_x \sub_\omega \mathcal{C}$ and ${\downarrow}x = [y, x] \cup {\downarrow}y \sub \bigcup \mathcal{V}_x$.
\end{proof}

The heart of the compactness proof is the following observation.

\begin{Proposition}\label{proposition : main proposition}
	For each ordinal $\alpha$ there exist $
		\mathcal{U}^\alpha \sub_\omega \mathcal{C}
	$ and an antichain
	$
		F^\alpha \sub_\omega X
	$
	such that
	\begin{equation}\label{eq 1 compactness}
		X \smallsetminus {\uparrow}F^\alpha \sub \bigcup \mathcal{U}^\alpha \text{ and there are no }
		x \in F^\alpha,
		\beta < \alpha,
		\text{ and }
		y \in F^\beta
		\text{ such that }
		x < y.
	\end{equation}
Furthermore, if $\alpha = \beta + 1$, then $F^{\alpha} \sub {\uparrow}F^\beta \smallsetminus F^\beta$.
\end{Proposition}

For the sake of readability, we will postpone the proof of the above proposition to the end of this section.  Instead, we shall now explain how the above proposition can be used to prove that $\mathcal{C}$ has a finite subcover.

\begin{Corollary}\label{Cor : very last cor : compactness : contradiction}
There exists an ordinal $\alpha$ such that for every ordinal $\gamma \geq \alpha$ it holds that $F^{\gamma} \sub \bigcup_{\beta < \gamma} F^\beta$.
\end{Corollary}

\begin{proof}
Suppose the contrary, i.e., that for every ordinal $\alpha$ there exists an ordinal $\alpha' \geq \alpha$ such that $F^{\alpha'} \centernot\sub \bigcup_{\beta < \alpha'} F^\beta$. Then for each ordinal $\alpha$ we define an ordinal $\alpha^*$ as follows. First, we let $0^* \coloneqq 0'$. Then consider an ordinal $\alpha > 0$ and assume that $\gamma^*$ has been defined for each $\gamma < \alpha$. We let
\[
\alpha^* \coloneqq \left(\sup\left(\{ \alpha \} \cup \{ \gamma^* : \gamma < \alpha\}\right)+1\right)'.
\]
It is easy to see that for every pair of ordinals $\alpha < \beta$ we have $\alpha^* < \beta^*$ and that for each ordinal $\alpha$,
\begin{equation}\label{Eq : Y alpha ast}
\alpha \leq \alpha^* \, \, \text{ and } \, \, F^{\alpha^*} \centernot\sub \bigcup_{\beta < \alpha^*} F^{\beta}.
\end{equation}
 
In view of the right hand side of the above display, for every ordinal $\alpha$ there exists 
\[
x_\alpha \in F^{\alpha^\ast} \smallsetminus \bigcup_{\beta < \alpha^\ast} F^\beta.
\] 
We will prove that $x_\alpha \ne x_\beta$ for each pair of distinct ordinals $\alpha$ and $\beta$. Suppose that $\alpha \ne \beta$. By symmetry we may assume $\alpha < \beta$. As we mentioned, this implies $\alpha^\ast < \beta^\ast$. As $x_\alpha \in F^{\alpha^\ast}$, we obtain $x_\alpha \notin F^{\beta^\ast} \smallsetminus \bigcup_{\gamma < \beta^\ast} F^\beta$. Since $x_\beta \in F^{\beta^\ast} \smallsetminus \bigcup_{\gamma < \beta^\ast} F^\beta$, we conclude that $x_\alpha \ne x_\beta$ as desired. Hence, $\{ x_\alpha : \alpha \text{ is an ordinal}\}$ is a proper class. But this contradicts the assumption that $X$ is a set containing each $x_\alpha$.
\end{proof}

We are now ready to show that $\mathcal{C}$ has a finite subcover. Suppose the contrary, with a view of contradiction, i.e., that 
\begin{equation}\label{eq : not compact}
	\text{ there is no } \mathcal{U} \sub_\omega \mathcal{C} \text{ such that } X \sub \bigcup \mathcal{U}.
\end{equation}

Recall from Corollary \ref{Cor : very last cor : compactness : contradiction} that there exists an ordinal $\alpha$ such that 
\[
F^{\gamma} \sub \bigcup_{\beta < \gamma} F^\beta\text{ for each ordinal }\gamma \geq \alpha.
\] 
We will show that the set $F^{\alpha+1}$ is nonempty. Suppose the contrary, with a view to contradiction. Then the left hand side of condition (\ref{eq 1 compactness}) yields $X \sub \bigcup \mathcal{U}^{\alpha+1}$.\ Therefore, the finite family $\mathcal{U} \coloneqq \mathcal{U}^{\alpha+1}$ contradicts condition (\ref{eq : not compact}). Hence, we conclude that  $F^{\alpha+1} \ne \emptyset$.

Then there exists $y \in F^{\alpha+1}$. In view of the above display, there also exists $\beta \leq \alpha$ such that $y \in F^\beta$. As $\alpha+1$ is a successor ordinal, the last part of Proposition \ref{proposition : main proposition} implies $y \in F^{\alpha+1} = {\uparrow}F^\alpha \smallsetminus F^\alpha$. Therefore, there exists $x \in F^\alpha$ such that $x < y$. Since $F^\beta$ is an antichain by Proposition \ref{proposition : main proposition} and $y \in F^\beta$, we obtain $x \notin F^\beta$. Together with $x \in F^\alpha$, this yields $\alpha \ne \beta$. Thus, from $\beta \leq \alpha$ it follows that $\beta < \alpha$. As $x \in F^\alpha$, $y \in F^\beta$, and $x < y$, this contradicts the right hand side of condition (\ref{eq 1 compactness}). Hence, we conclude that $\mathcal{C}$ has a finite subcover as desired. Therefore, in order to establish Theorem \ref{Thm : Esakia compactness : trees}, it only remains to prove Proposition \ref{proposition : main proposition}.

\begin{proof}[Proof of Proposition \ref{proposition : main proposition}.]
As $\mathcal{C}$ covers $X$, for each $x \in X$ there exists $U_x \in \mathcal{C}$ such that $f_x\left(\h{X}\right) \in U_x$. Since $\mathcal{C} \subseteq \mathcal{S}_{\h{X}}$, we can apply the Main Lemma \ref{lemma : fundamental lemma}, obtaining $v_x \in X$ such that $v_x \leq x$ and  $Y_x \sub_\omega {\uparrow}v_x \cap X_{>\h{x}}$ and $Z_x \sub_\omega X_{<\h{X}} \cap {\uparrow} v_x$ such that
\begin{equation}\label{Eq : compactness : main lemma property}
{\uparrow}v_x \smallsetminus \left({\uparrow}Y_x \cup {\downarrow}Z_x\right) \sub U_x.
\end{equation}
Furthermore, $\h{v_x}$ is either zero or a successor ordinal. In addition, for each $x \in X$ there exists $\mathcal{V}_x \sub_\omega \mathcal{C}$ such that 
\begin{equation}\label{Eq : compactness : Vx property}
{\downarrow}x \subseteq \bigcup \mathcal{V}_x
\end{equation} 
by Proposition \ref{proposition : covering X minus down Z suffices}. The objects $v_x, U_x, Y_x, Z_x$, and $\mathcal{V}_x$ will be used repeatedly in the proof, which proceeds by induction on $\alpha$.

\subsection*{Base case.}
If $\alpha = 0$, we let $\mathcal{U}^0 \coloneqq \emptyset$ and define $F^0$ as the singleton containing the root of $X$. Then $X \smallsetminus {\uparrow}F^0 = X \smallsetminus X =  \emptyset \sub \bigcup \mathcal{U}^0$ and the other conditions in the statement of Proposition \ref{proposition : main proposition} are clearly satisfied.

\subsection*{Successor case.} Consider a successor ordinal $\alpha + 1$. By the inductive hypothesis there exist $
		\mathcal{U}^\alpha \sub_\omega \mathcal{C}
	$ and an antichain
	$
		F^\alpha \sub_\omega X
	$
satisfying condition (\ref{eq 1 compactness}). We let 	
\[
\begin{array}{rcl}
A^{\alpha+1}& \coloneqq &
			\left\{ x : x \in Y_y \cap {\uparrow}y \text{ for some }y \in F^\alpha \text{ and there are no }\beta \leq \alpha \text{ and } z \in F^\beta \text{ s.t. }x < z\right\};\\
	F^{\alpha+1}
		& \coloneqq &
			\min A^{\alpha+1};\\
	\mathcal{U}^{\alpha+1}
		& \coloneqq &
			\mathcal{U}^\alpha \, \cup \, \left\{U_x \colon x \in  F^\alpha\right\} \, \cup \, \left\{U : \text{there exist }y \in F^\alpha \text{ and }x \in Z_y \text{ s.t. }U \in \mathcal{V}_x  \right\}.
\end{array}
\]

As $F^\alpha$ is finite and so is $Y_y$ for each $y \in F^\alpha$, the set $A^{\alpha+1}$ is also finite. Consequently, $F^{\alpha+1}$ is a finite antichain. On the other hand, as $\mathcal{U}^\alpha$ and $F^\alpha$ are finite and so is $Z_y$ for each $y \in F^\alpha$ as well as $\mathcal{V}_x$ for each $x \in Z_y$, the set $\mathcal{U}^{\alpha+1}$ is also finite. Furthermore, $\mathcal{U}^{\alpha+1} \subseteq \mathcal{C}$ because $\mathcal{U}^\alpha \subseteq \mathcal{C}$ by the inductive hypothesis and $\{U_x\} \cup \mathcal{V}_x \sub \mathcal{C}$ for each $x \in X$ by assumption. Hence, $\mathcal{U}^{\alpha+1}\sub_\omega \mathcal{C}$ and $F^{\alpha+1} \sub_\omega X$, where $F^{\alpha+1}$ is also an antichain. Therefore, it only remains to prove that $\mathcal{U}^{\alpha+1}$ and $F^{\alpha+1}$ satisfy condition  (\ref{eq 1 compactness}) and the last part of Proposition \ref{proposition : main proposition}.

\begin{Claim}\label{Claim : compactness : claim number 1}
We have $X \smallsetminus {\uparrow} F^{\alpha+1} \sub \bigcup \mathcal{U}^{\alpha+1}$.
\end{Claim}
\begin{proof}[Proof of the Claim]
Let $x \in X \smallsetminus {\uparrow} F^{\alpha+1}$. By the inductive hypothesis we have $X \smallsetminus {\uparrow} F^{\alpha} \sub \bigcup \mathcal{U}^{\alpha}$. Therefore, if $x \notin {\uparrow} F^{\alpha}$, then $x \in X \smallsetminus {\uparrow} F^{\alpha} \sub \bigcup \mathcal{U}^{\alpha} \subseteq \bigcup \mathcal{U}^{\alpha+1}$, where the last inclusion follows from the assumption that $\mathcal{U}^\alpha \subseteq \mathcal{U}^{\alpha+1}$. Then we consider the case where $x \in {\uparrow} F^{\alpha}$. There exists $y \in F^\alpha$ such that $y \leq x$.\ We have two cases:\ either $x \in {\downarrow}Z_y$ or $x \notin {\downarrow}Z_y$. First, suppose that $x \in {\downarrow}Z_y$. Then there exists $z \in Z_y$ such that $x \leq z$. By condition (\ref{Eq : compactness : Vx property}) we have ${\downarrow} z \subseteq \bigcup \mathcal{V}_z$. Since $x \leq z$, this yields $x \in\bigcup \mathcal{V}_z$. On the other hand, from $z \in Z_y$ and $y \in F^\alpha$ it follows that $\mathcal{V}_z \sub \mathcal{U}^{\alpha+1}$. Hence, $x \in \bigcup \mathcal{U}^{\alpha+1}$ as desired. Then we consider the case where $x \notin {\downarrow}Z_y$. Again, we have two cases: either $x \notin {\uparrow}Y_y$ or $x \in {\uparrow}Y_y$. First, suppose that $x \notin {\uparrow}Y_y$. Together with $v_y \leq y \leq x$ and $x \notin {\downarrow}Z_y$, this yields $x \in {\uparrow}v_y \smallsetminus \left({\uparrow}Y_y \cup {\downarrow}Z_y\right)$. By condition (\ref{Eq : compactness : main lemma property}) this implies $x \in U_y$. Since $y \in F^\alpha$, we have $U_y \in \mathcal{U}^{\alpha+1}$ and, therefore, $x \in \bigcup  \mathcal{U}^{\alpha+1}$ as desired. It only remains to consider the case where $x \in {\uparrow}Y_y$. We will show that this cases never happens, in the sense that it leads to a contradiction. First, as $x \in {\uparrow}Y_y$, there exists $z \in Y_y$ such that $z \leq x$. We will prove that $z \in A^{\alpha+1}$. Since $X$ is a tree and $y, z \leq x$, the elements $y$ and $z$ must be comparable. As $Y_y \sub X_{>\h{y}}$ and $z \in Y_y$, we deduce $y < z$. Therefore, $z \in Y_y \cap {\uparrow}y$ and $y \in F^\alpha$. Consequently, to prove that $z \in A^{\alpha+1}$, it only remains to show that there are no $\beta \leq \alpha$ and $w \in F^\beta$ such that $z < w$. Suppose, on the contrary, that there exist such $\beta$ and $w$. From $y < z < w$ it follows that $y < w$. Recall that $\beta \leq \alpha$. Then either $\beta < \alpha$ or $\beta = \alpha$. The case where $\beta < \alpha$ cannot happen because $F^\alpha$ satisfies the right hand side of condition (\ref{eq 1 compactness}) and $y \in F^\alpha, w \in F^\beta$, and $y <  w$. Therefore, we obtain $\alpha = \beta$. As a consequence, $y, w \in F^\alpha$ because $y \in F^\alpha$ and $w \in F^\beta$. Together with $y < w$, this contradicts the assumption that $F^\alpha$ is an antichain. Hence, we conclude that $z \in A^{\alpha+1}$. Since the set $A^{\alpha+1}$ is finite and $F^{\alpha+1} = \min A^{\alpha+1}$, this yields $z \in {\uparrow}F^{\alpha+1}$. As $z \leq x$, we obtain $x \in {\uparrow}F^{\alpha+1}$, a contradiction with the assumption that $x \in X \smallsetminus {\uparrow}F^{\alpha+1}$.
\end{proof}

By the Claim \ref{Claim : compactness : claim number 1} the set $F^{\alpha+1}$ satisfies the left hand side of condition (\ref{eq 1 compactness}). The right hand side of the same conditions holds by the  definition of $F^{\alpha+1}$. Therefore, it only remains to prove the last part of Proposition \ref{proposition : main proposition}, namely, that $F^{\alpha+1} \sub {\uparrow} F^\alpha \smallsetminus F^\alpha$. To this end, consider $x \in F^{\alpha+1}$. By the definition of $F^{\alpha+1}$ we have $x \in Y_y \cap {\uparrow}y$ for some $y \in F^\alpha$. Therefore, $x \in {\uparrow}F^\alpha$. It only remains to prove that $x \notin F^\alpha$. From $x \in Y_y \sub X_{> \h{y}}$ and $x \geq y$ it follows that $x < y$. As $F^\alpha$ is an antichain containing $y$, this implies $x \notin F^\alpha$.

\subsection*{Limit case.} For each nonempty $Y \subseteq X$ let
\[
\sup{\!^\ast} Y \coloneqq \{ x \in X : x \text{ is the supremum of a maximal chain }Z \subseteq Y \}.
\]

Suppose that $\alpha$ is a limit ordinal. By the inductive hypothesis for each $\beta < \alpha$ there exist $
		\mathcal{U}^\beta \sub_\omega \mathcal{C}
	$ and an antichain
	$
		F^\beta \sub_\omega X
	$
satisfying condition (\ref{eq 1 compactness}). We let
\[
F \coloneqq \bigcup_{\beta < \alpha} F^\beta \cup \bigcup_{\beta < \alpha}\left(X \smallsetminus {\uparrow}F^\beta\right) \, \, \text{ and } \, \, F^\ast \coloneqq {\downarrow}\left(\sup{\!^{\ast}} F\right).
\]
Notice that $F$ is nonempty because it contains the root of $X$ (the latter belongs to $F^0$ by construction and $F^0 \subseteq F$). Therefore, every maximal chain in $F$ is nonempty and, therefore, has a supremum in $X$ by assumption.

The proof relies on a series of technical observation.

\begin{Claim}\label{Claim : downsets F and F ast}
The sets $F$ and $F^*$ are nonempty downsets of $X$.
\end{Claim}

\begin{proof}[Proof of the Claim]
We begin by proving that $F$ is a nonempty downset of $X$. First, $F$ is nonempty because $0 < \alpha$ and $F^0 \subseteq F$ is the singleton containing the root of $X$. To prove that $F$ is a downset, for every ordinal $\beta < \alpha$ let
\[
G^\beta \coloneqq F^\beta \cup \left(X \smallsetminus {\uparrow}F^\beta\right).
\]
We show that each $G^\beta$ is a downset. Consider $x \in G^\beta$ and $y < x$. We need to prove that $y \in G^\beta$. There are two cases: either $x \in F^\beta$ or $x \in X \smallsetminus {\uparrow}F^\beta$. Suppose that $x \in F^\beta$. Then $y$ cannot belong to ${\uparrow}F^\beta$, otherwise there exists $z \in F^\beta$ such that $z \leq y < x$. Since $x, z \in F^\beta$, this contradicts the assumption that $F^\beta$ is an antichain. Hence, $y \in X \smallsetminus {\uparrow}F^\beta \subseteq G^\beta$ as desired. Then we consider the case where $x \in X \smallsetminus {\uparrow}F^\beta$. Since $X \smallsetminus {\uparrow}F^\beta$ is a downset and $y \leq x$, we obtain $y \in X \smallsetminus {\uparrow}F^\beta \subseteq G^\beta$ too. Hence, each $G^\beta$ is a downset. As $F = \bigcup_{\beta < \alpha}G^\beta$, we conclude that $F$ is a downset. Lastly, $F^*$ is a nonempty downset by definition.
\end{proof}

\begin{Claim}\label{Claim : F ast is a good poset}
The poset $\langle F^*; \leq \rangle$ is order compact and each of its nonempty chains has a supremum in $F^*$.
\end{Claim}

\begin{proof}[Proof of the Claim]
Recall that $F^*$ is a nonempty downset of $X$ by Claim \ref{Claim : downsets F and F ast}. Together with the assumption that $X$ is a tree with enough gaps, this yields that $F^*$ is also a tree with enough gaps. We will show that each of its nonempty chains has a supremum in $F^*$.\ Together with Theorem \ref{Thm:Lewis-tree}, this  implies that $X$ is representable and, therefore, order compact by Proposition \ref{Prop:properties}. As such, in order to conclude the proof, it suffices to show that in $F^*$ each nonempty chain has a supremum.

All suprema in the rest of the proof will be computed in $X$, unless said otherwise.\ Consider a nonempty chain $C = \{ x_i : i \in I \}$ in $F^*$. We will prove that $C$ has a supremum in the poset $\langle F^*; \leq \rangle$. We begin by showing that
\begin{equation}\label{Eq : zi is the sup of F}
x_i = \sup \left(F \cap {\downarrow}x_i\right) \, \, \text{ for each } i \in I.
\end{equation}
Consider $i \in I$. If $x_i \in F$, clearly $x_i = \sup \left(F \cap {\downarrow}x_i\right)$ and we are done. Then consider the case where $x_i \notin F$. Since $x_i \in F^* = {\downarrow}\left(\sup^*F\right)$, there exists $y \in \sup^* F$ such that $x_i \leq y$. Furthermore, $y = \sup \left(F \cap {\downarrow}y\right)$ because $y \in \sup^* F$. We have two cases: either $x_i = y$ or $x_i \ne y$. If $x_i = y$, we have $x_i = y = \sup \left(F \cap {\downarrow}y\right) = \sup \left(F \cap {\downarrow}x_i\right)$ and we are done. Then we consider the case where $x_i \ne y$. As $x_i \leq y$, we have $x_i < y$. Since $y = \sup \left(F \cap {\downarrow}y\right)$, this guarantees the existence of $z \in F$ such that $z \leq y$ and $z \nleq x_i$. As $X$ is a tree and $x_i, z \leq y$, the elements $x_i$ and $z$ must be comparable. Together with $z \nleq x_i$, this yields $x_i \leq z$. Since $z \in F$ and $F$ is a downset by Claim \ref{Claim : downsets F and F ast}, we obtain $x_i \in F$, a contradiction. This establishes condition (\ref{Eq : zi is the sup of F}).

 Then consider the set
\[
D \coloneqq \bigcup_{i \in I}\left(F \cap {\downarrow}x_i\right).
\]
We will prove that $D$ is a nonempty chain. First, recall that the root of $X$ belongs to $F^0$ and, therefore, to $F$ by construction. Thus, $D$ is nonempty. Then consider $y, z \in D$. By the definition of $D$ there exist $i, j \in I$ such that $y \leq x_i$ and $z \leq x_j$. Since $C$ is a chain, by symmetry we may assume that $x_i \leq x_j$. Therefore, $y, z \leq x_j$. As $X$ is a tree, we conclude that $y$ and $z$ are comparable. Hence, $D$ is a nonempty chain as desired. Consequently, $\sup D$ exists by the assumptions on $X$. Together with the definitions of $C$ and $D$ and with condition (\ref{Eq : zi is the sup of F}), this yields that also $\sup C$ exists and coincides with $\sup D$. Furthermore, the definition of $D$ ensures that $D \sub F$. Since $D$ can be extended to a maximal chain in $F$ by Zorn's Lemma, we obtain $\sup C = \sup D \in {\downarrow}\left( \sup^* F\right) = F^*$. Thus, the supremum of $C$ computed in $X$ exists and belongs to $F^*$. Clearly, this coincides with the supremum of $C$ computed in $F^*$. Therefore, we conclude that $C$ has a supremum also in the poset $\langle F^*; \leq \rangle$ as desired.
\end{proof}

Recall that for each $x \in X$ we have $v_x \leq x$ and that $\h{x}$ is either zero or a successor ordinal.

\begin{Claim}\label{Claim : sup Y}
	For every $x \in \sup^* F$ there exist an ordinal $\gamma_x < \alpha$, an element $y_x \in X$, an order open subset $V_x$ of $\langle F^*; \leq \rangle$, and $\mathcal{W}_x \sub_\omega \mathcal{C}$ satisfying the following conditions:
	\benroman
	\item\label{item : 1 : the new items : comp X} $v_x \leq y_x \leq x$;
	\item\label{item : 2 : the new items : comp X} $V_x$ is disjoint both from ${\uparrow}\left(F^{\gamma_x} \smallsetminus {\uparrow}y_x\right)$ and ${\uparrow}Y_x$;
	\item\label{item : 3 : the new items : comp X} $x \in V_x \sub U_x \cup \bigcup \mathcal{W}_x \cup \bigcup \mathcal{U}^{\gamma_x}$.
	\eroman
	\end{Claim}

\begin{proof}[Proof of the Claim]
	Consider $x \in \sup^\ast F$. We will prove that there exist
	\begin{equation}\label{Eq : an easy to forget condition}
\gamma_x < \alpha \, \, \text{ and } \, \, y_x \in F^{\gamma_x} \cup \left(X \smallsetminus {\uparrow}F^{\gamma_x}\right) \, \, \text{ such that }\,\,	v_x \leq y_x \leq x.
	\end{equation}
This will establish condition (\ref{item : 1 : the new items : comp X}). Recall that $x$ is the supremum of a maximal chain of $\langle F; \leq \rangle$ because $x \in \sup^\ast F$. We have two cases: either $\h{x}$ is a limit ordinal or not. First, suppose that $\h{x}$ is not a limit ordinal. Since $x$ is the supremum of a nonempty chain of $\langle F; \leq \rangle$, this yields $x \in F$. Consequently, there exists $\gamma_x < \alpha$ such that $x \in F^{\gamma_x} \cup \left(X \smallsetminus {\uparrow}F^{\gamma_x}\right)$. Therefore, letting $y_x \coloneqq x$, we are done. Then we consider the case where $\h{x}$ is a limit ordinal. As $\h{v_x}$ is either zero or a successor ordinal and $v_x \leq x$, this yields $v_x < x$. Since $x$ is the supremum of chain of $\langle F; \leq \rangle$, there exist $\gamma_x < \alpha$ and $y_x \in F^{\gamma_x}\cup \left(X \smallsetminus {\uparrow}F^{\gamma_x}\right)$ such that $y_x \leq x$ and $y_x \nleq v_x$. As $X$ is a tree and $y_x, v_x \leq x$, we deduce that either $y_x \leq v_x$ or $v_x \leq y_x$. Together with $y_x \nleq v_x$, this yields $v_x \leq y_x$ and establishes the above display.

Let 
	\[
V_x' \coloneqq		{\uparrow}
			\left(
				F^{\gamma_x} \smallsetminus {\uparrow}y_x
			\right)
		\cup
			{\uparrow}Y_x
		\cup 
			{\downarrow}
				\left(
					Z_x \smallsetminus{\uparrow}x
				\right) \,\, \text{ and } \, \, V_x \coloneqq F^* \smallsetminus V_x'.
	\]
Notice that $V_x$ satisfies condition (\ref{item : 2 : the new items : comp X}) by definition. We will prove that $V_x$ is an order open set of the poset $\langle F^*; \leq \rangle$. First, observe that $V_x \subseteq F^*$ by the definition of $V_x$. Then consider the sets
\[
A \coloneqq F^* \cap \left(\left(F^{\gamma_x} \smallsetminus {\uparrow}y_x\right) \cup Y_x\right) \, \, \text{ and } \, \, B \coloneqq \max \left(F^* \cap {\downarrow}\left(Z_x \smallsetminus {\uparrow}z\right)\right).
\]
We shall see that $A, B \sub_\omega F^*$. Since $F^{\gamma_x}$ and $Y_x$ are fine, we obtain $A \sub_\omega F^*$. On the other hand, as $X$ is a tree and $Z_x$ finite, the set ${\downarrow}\left(Z_x \smallsetminus {\uparrow}z\right)$ is the union of $n$ chains for some nonnegative integer $n$.\ Consequently, $\vert B\vert \leq n$ and, therefore, $B \sub_\omega F^*$ as desired. From $A, B \sub_\omega F^*$ and Lemma \ref{Lem : typical order open sets} it follows that $F^* \smallsetminus \left({\uparrow}A \cup {\downarrow}B\right)$ is an order open set of $\langle F^*; \leq \rangle$. 

To prove that $V_x$ is also an order open set of $\langle F^*; \leq \rangle$, we rely on the equalities
\begin{equation}\label{Eq : the set B has enough maximal elements}
F^* \cap \left({\uparrow}\left(
				F^{\gamma_x} \smallsetminus {\uparrow}y_x
			\right)\cup
			{\uparrow}Y_x\right) =F^* \cap  {\uparrow}A \, \, \text{ and }\, \, F^* \cap {\downarrow}\left(Z_x \smallsetminus{\uparrow}x\right) = F^* \cap {\downarrow}B.
\end{equation}
First, observe that
\begin{align*}
F^* \cap \left({\uparrow}\left(
				F^{\gamma_x} \smallsetminus {\uparrow}y_x
			\right)
		\cup
			{\uparrow}Y_x\right)
		 &= F^* \cap {\uparrow}\left(
				\left(F^{\gamma_x} \smallsetminus {\uparrow}y_x\right) \cup Y_x
			\right)\\
		&= F^* \cap  {\uparrow}\left(F^* \cap \left(\left(F^{\gamma_x} \smallsetminus {\uparrow}y_x\right) \cup Y_x\right)\right) \\
		&=F^* \cap  {\uparrow}A,
\end{align*}
where the first equality is straightforward, the second holds because $F^*$ is a downset of $X$, and the third holds by the definition of $A$. This establishes the left hand side of condition (\ref{Eq : the set B has enough maximal elements}). Then we turn to prove the right hand side of the same condition. The inclusion from right to left is an immediate consequence of the definition of $B$. To prove the other inclusion, consider $z \in F^* \cap {\downarrow}\left(Z_x \smallsetminus{\uparrow}x\right)$. By Zorn's lemma there exists a maximal chain $C \sub F^* \cap {\downarrow}\left(Z_x \smallsetminus{\uparrow}x\right)$ such that $z \in C$. Since $C$ is a nonempty chain of $F^*$, it has a supremum $\sup C$ in $\langle F^*; \leq \rangle$ by Claim \ref{Claim : F ast is a good poset}. We will prove that $\sup C \in F^*\cap {\downarrow}\left(Z_x \smallsetminus{\uparrow}x\right)$. Since $\sup C \in F^*$, it suffices to show that $\sup C \in F^* \in {\downarrow}\left(Z_x \smallsetminus{\uparrow}x\right)$. Recall that $Z_x$ is finite. Therefore, so is $Z_x \smallsetminus{\uparrow}x$. Furthermore, $Z_x \smallsetminus{\uparrow}x$ is nonempty because $z \in {\downarrow}\left(Z_x \smallsetminus{\uparrow}x\right)$. Then consider an enumeration $Z_x \smallsetminus{\uparrow}x = \{ z_1, \dots, z_n \}$. We will show that $C \subseteq {\downarrow}z_i$ for some $i \leq n$. Suppose the contrary, with a view to contradiction. Then for each $i \leq n$ there exists $c_i \in C$ such that $c_i \nleq z_i$. As $C$ is a chain, the set $\{ c_1, \dots, c_n \}$ has a maximum $c$. Clearly, we have $c \nleq z_1, \dots, z_n$, a contradiction with the assumption that $C \sub {\downarrow}\left(Z_x \smallsetminus{\uparrow}x\right)$. Hence, there exists $i \leq n$  such that $C \subseteq {\downarrow}z_i$. Consequently, $\sup C \leq z_i$. Since $z_i \in Z_x \smallsetminus{\uparrow}x$, we obtain $\sup C \in {\downarrow}\left(Z_x \smallsetminus{\uparrow}x\right)$ as desired. From $\sup C \in \left(F^* \cap {\downarrow}\left(Z_x \smallsetminus{\uparrow}x\right)\right)$ and the maximality of the chain $C$ it follows that $\sup C \in \max \left(F^* \cap {\downarrow}\left(Z_x \smallsetminus{\uparrow}x\right)\right) = B$. Together with $z \in C$, this yields $z \in {\downarrow}B$. As $z \in F^*$, we conclude that $z \in F^* \cap {\downarrow}B$, establishing condition (\ref{Eq : the set B has enough maximal elements}).

Lastly, observe that
\begin{align*}
V_x &= F^* \smallsetminus \left({\uparrow}
			\left(
				F^{\gamma_x} \smallsetminus {\uparrow}y_x
			\right)
		\cup
			{\uparrow}Y_x
		\cup 
			{\downarrow}
				\left(
					Z_x \smallsetminus{\uparrow}x
				\right)\right)\\
				&=F^* \smallsetminus \left(\left(F^* \cap {\uparrow}\left(\left(
				F^{\gamma_x} \smallsetminus {\uparrow}y_x\right)
		\cup
			Y_x\right)\right)
		\cup \left(F^* \cap
			{\downarrow}
				\left(
					Z_x \smallsetminus{\uparrow}x
				\right)\right)\right)\\
				&= F^* \smallsetminus \left(\left(F^* \cap {\uparrow}A\right) \cup \left(F^* \cap {\downarrow}B\right)\right)\\
				&= F^* \smallsetminus \left({\uparrow}A \cup {\downarrow}B\right),
\end{align*}
where the first equality holds by the definition of $V_x$, the second and the last are straightforward, and the third holds by condition (\ref{Eq : the set B has enough maximal elements}). Therefore, since $F^* \smallsetminus \left({\uparrow}A \cup {\downarrow}B\right)$ is an order open set of $\langle F^*; \leq \rangle$, we conclude that so is $V_x$.

Therefore, it only remains to construct $\mathcal{W}_x \sub_\omega \mathcal{C}$ so that condition (\ref{item : 3 : the new items : comp X}) holds. Let
\[
\mathcal{W}_x \coloneqq \{ U : U \in \mathcal{V}_z \text{ for some }z \in Z_x \}.
\]
Since $Z_x$ is finite and $\mathcal{V}_z \sub_\omega \mathcal{C}$ for each $z \in Z_x$, we obtain $\mathcal{W}_x \sub_\omega \mathcal{C}$. Then we turn to prove condition (\ref{item : 3 : the new items : comp X}). 

We begin by showing that $x \in V_x$. Suppose the contrary, with a view to contradiction. Since $x \in \sup^\ast F \subseteq F^*$ by assumption, we obtain $x \in F^* \smallsetminus V_x \sub V_x'$. From the definition of $V_x'$ it follows that
\[
\text{either }\, \, x \in {\uparrow}
			\left(
				F^{\gamma_x} \smallsetminus {\uparrow}y_x 
			\right) \, \, \text{ or } \, \, x \in 
			{\uparrow}Y_x \, \, \text{ or } \, \, 
	x \in		{\downarrow}
				\left(
					Z_x \smallsetminus{\uparrow}x
				\right).
\]
	First, suppose $x \in {\uparrow}\left(F^{\gamma_x} \smallsetminus {\uparrow}y_x \right)$. Then there exists $z \in F^{\gamma_x} \smallsetminus {\uparrow}y_x$ such that $z \leq x$. Since $X$ is a tree and $y_x, z \leq x$ (for $y_x \leq x$, see condition (\ref{Eq : an easy to forget condition})), we deduce that either $z \leq y_x$ or $y_x \leq z$.  As $z \in F^{\gamma_x} \smallsetminus {\uparrow}y_x$, this amounts to $z < y_x$. In view of condition (\ref{Eq : an easy to forget condition}), either $y_x \in F^{\gamma_x}$ or $y_x \in X \smallsetminus {\uparrow}F^{\gamma_x}$. We will show that both cases lead to a contradiction. If  $y_x \in F^{\gamma_x}$, we have $y_x, z \in F^{\gamma_x}$. Together with $z < y_x$, this contradicts the assumption that $F^{\gamma_x}$ is an antichain. On the other hand, if $y_x \in X \smallsetminus {\uparrow}F^{\gamma_x}$, we obtain a contradiction with $z < y_x$ and $z \in F^{\gamma_x}$. Lastly, the case where $x \in {\uparrow} Y_x$ leads to a contradiction because $Y_x \sub X_{>\h{x}}$, and the case  
	$x \in {\downarrow}
				\left(
					Z_x \smallsetminus{\uparrow}x
				\right)$ is obviously impossible. Hence, we conclude that $x \in V_x$.

 Therefore, to conclude the proof, it only remains to show that
\[
V_x \sub U_x \cup \bigcup \mathcal{W}_x \cup \bigcup \mathcal{U}^{\gamma_x}.
\]
Consider $y \in V_x$. There are two cases: either $y \in  {\downarrow} Z_x$ or $y \notin {\downarrow}Z_x$. First, suppose that $y \in  {\downarrow} Z_x$. Then there exists $z \in Z_x$ such that $y \leq z$. Therefore, $\mathcal{V}_z \subseteq \mathcal{W}_x$ by the definition of $\mathcal{W}_x$. From condition (\ref{Eq : compactness : Vx property}) and $y \leq z$ it follows that $y \in {\downarrow}z \subseteq \bigcup \mathcal{V}_z \subseteq \bigcup \mathcal{W}_x$ as desired. Then we consider the case where $y \notin {\downarrow}Z_x$. Again, we have two cases: either $y \notin {\uparrow} F^{\gamma_x}$ or $y \in {\uparrow} F^{\gamma_x}$. If $y \notin {\uparrow} F^{\gamma_x}$, we have $y \in X \smallsetminus {\uparrow} F^{\gamma_x}$. Therefore, the fact that $\mathcal{U}^{\gamma_x}$ and $F^{\gamma_x}$ satisfy condition (\ref{eq 1 compactness}) ensures that $y \in \bigcup \mathcal{U}^{\gamma_x}$ and we are done. 
	Lastly, we consider the case where $y \in {\uparrow} F^{\gamma_x}$. Since $y \in V_x$ by assumption and $V_x \subseteq \left({\uparrow} Y_x\right)^c \cap \left({\uparrow}\left(F^{\gamma_x} \smallsetminus {\uparrow} y_x\right)\right)^c$ by the definition of $V_x$, we have $y \notin {\uparrow} Y_x$ and $y \notin {\uparrow}\left(F^{\gamma_x} \smallsetminus {\uparrow} y_x\right)$. Together with $y \in {\uparrow} F^{\gamma_x}$, the latter yields $y \in {\uparrow}y_x$. Therefore, $y \in {\uparrow}y_x$, $y \notin {\uparrow} Y_x$, and $y \notin {\downarrow}Z_x$. Since $v_x \leq y_x$ by condition (\ref{Eq : an easy to forget condition}), this yields	$y \in {\uparrow}v_x \smallsetminus \left({\uparrow}Y_x \cup {\downarrow}Z_x\right)$. By condition (\ref{Eq : compactness : main lemma property}) we conclude that $y \in U_x$ as desired.
\end{proof}
			
Recall that Claim \ref{Claim : sup Y} associates a set $V_x$ with every $x \in \sup^* F$. Using these sets, we obtain the following:
	
\begin{Claim}\label{Claim : cantor compactness}
	There exist $G \sub_\omega \sup^\ast F$ and $\Gamma \sub_\omega \alpha$ such that 
	\[
	F^*
			\sub
		\bigcup\limits_{x \in G}
			V_{x}
			\cup
		\bigcup\limits_{\beta \in \Gamma}\left(
			X \smallsetminus {\uparrow} F^{\beta}\right).
	\]
\end{Claim}
\begin{proof}[Proof of the Claim.]
	First, we show that 
	\begin{equation}\label{Eq : the application of order compactness in the c proof}
		F^*
			\sub
		\bigcup_{x \in \sup^\ast F}
			V_{x}
			\cup
		\bigcup_{\beta < \alpha}
			\left(X \smallsetminus {\uparrow} F^{\beta}\right).
	\end{equation}
To prove this, consider $y \in F^* = {\downarrow}\left(\sup^* F\right)$.\ Then there exists $x \in \sup^\ast F$ such that $y \leq x$. If $y = x$, Claim \ref{Claim : sup Y}(\ref{item : 3 : the new items : comp X}) ensures $y \in V_y$ and we are done. Then we consider the case where $y < x$. Since $x$ is the supremum of a maximal chain of $\langle F; \leq \rangle$, there exist $\beta < \alpha$ and $z \in F^\beta \cup \left(X \smallsetminus {\uparrow}F^{\beta}\right)$ such that $z \leq x$ and $z \nleq y$. As $X$ is a tree and $y, z \leq x$, the elements $y$ and $z$ must be comparable. Together with $z \nleq y$, this yields $y < z$. We will prove that $y \notin {\uparrow}F^{\beta}$. Recall that $z \in F^\beta \cup \left(X \smallsetminus {\uparrow}F^{\beta}\right)$. We will consider the cases where 
$z \in F^\beta$ and $z \in X \smallsetminus {\uparrow}F^{\beta}$ separately. First, suppose that $z \in F^\beta$. Since $F^\beta$ is an antichain containing $z$ and $y < z$, we obtain $y \notin {\uparrow}F^\beta$ as desired. On the other hand, if $z \in X \smallsetminus {\uparrow}F^\beta$, then $y \notin {\uparrow}F^\beta$ because $y < z$. This concludes the proof that $y \notin {\uparrow}F^\beta$. Therefore, $y \in X \smallsetminus {\uparrow}F^\beta$ with $\beta < \alpha$, establishing the above display.

Now, observe that the following are order open sets of $\langle F^*; \leq \rangle$:
\benroman
\item\label{item : order opens : compactness : set 1} $V_x$ for each $x \in \sup^\ast F$;
\item\label{item : order opens : compactness : set 2} $F^* \cap \left(X \smallsetminus {\uparrow}F^\beta\right)$ for each $\beta < \alpha$.
\eroman
The sets in condition (\ref{item : order opens : compactness : set 1}) are order open by Claim \ref{Claim : sup Y}. To prove that the sets in condition (\ref{item : order opens : compactness : set 2}) are also order open, consider $\beta < \alpha$. Since $F^*$ is a downset of $X$, we have $F^* \cap \left(X \smallsetminus {\uparrow}F^\beta\right) = \left({\uparrow}\left(F^\beta \cap F^*\right)\right)^c$, where upsets and complements are computed in $F^*$. Therefore, it suffices to show that $\left({\uparrow}\left(F^\beta \cap F^*\right)\right)^c$ is an order open set of $\langle F^*; \leq \rangle$. The latter follows from Lemma \ref{Lem : typical order open sets} and the fact that $F^\beta \cap F^*$ is finite (because so is $F^\beta$).

Since the sets in conditions (\ref{item : order opens : compactness : set 1}) and (\ref{item : order opens : compactness : set 2}) are order open sets of $\langle F^*; \leq \rangle$ and this poset is order compact by Claim \ref{Claim : F ast is a good poset}, from condition (\ref{Eq : the application of order compactness in the c proof}) it follows that there exist $G \sub_\omega \sup^\ast F$ and $\Gamma \sub_\omega \alpha$ satisfying the statement of the claim.
\end{proof}

Using the sets $G$ and $\Gamma$ given by Claim \ref{Claim : cantor compactness} and the  sets $\mathcal{W}_x$ and the ordinals $\gamma_x$ given by Claim \ref{Claim : sup Y}, we let
	\begin{align*}
	A^\alpha &\coloneqq \{ x : x \in Y_y \text{ for some }y \in G \text{ and there are no }\beta < \alpha \text{ and }
									z \in F^\beta
								\text{ s.t. }
									x < z
					\};\\
				F^\alpha
			&\coloneqq
				\min A^\alpha;\\
				\mathcal{U}^\alpha & \coloneqq 
				\left\{ U : \text{there are }y \in G \text{ and }x \in Z_y \text{ s.t. }U \in \mathcal{V}_x \right\}	\cup \left\{ U : U \in \mathcal{W}_x \text{ for some }x   \in G\right\} \cup\\
				& \quad \, \,
					\left\{ U : U \in \mathcal{U}^{\gamma_x} \text{ for some }x   \in G \right\}
					\cup
				\left\{ U : U \in \mathcal{U}^{\beta} \text{ for some } \beta \in \Gamma\right\} \cup \left\{U_x \colon x \in G \right\}.
	\end{align*}
Since $G$ is finite by Claim \ref{Claim : cantor compactness} and so is $Y_y$ for each $y \in G$, the set $A^\alpha$ is also finite. Consequently, $F^\alpha$ is a finite antichain.
Moreover, $\mathcal{U}^\alpha$ is finite because so are the sets of the form $Z_x$,  $\mathcal{V}_x$, $\mathcal{W}_x$, and $\mathcal{U}^\beta$ for each $x \in X$ and $\beta < \alpha$ (for the case of $\mathcal{W}_x$, see Claim \ref{Claim : sup Y}) as well as the sets $G$ and $\Gamma$ by Claim \ref{Claim : cantor compactness}. Furthermore, $\mathcal{U}^\alpha \sub \mathcal{C}$ because $\mathcal{V}_x, \mathcal{W}_x, \mathcal{U}^\beta, \{U_x\} \subseteq \mathcal{C}$ for each $x \in X$ and $\beta < \alpha$ (for the case of $\mathcal{W}_x$, see Claim \ref{Claim : sup Y}). Hence, $\mathcal{U}^{\alpha}\sub_\omega \mathcal{C}$ and $F^{\alpha} \sub_\omega X$, where $F^{\alpha}$ is also an antichain. 

Observe that the last part of Proposition \ref{proposition : main proposition} holds vacuously because $\alpha$ is a limit ordinal. Therefore, it only remains to prove condition (\ref{eq 1 compactness}). The right hand side of this condition holds by the definition of $A^\alpha$ and the fact that $F^\alpha \sub A^\alpha$. Therefore, we turn to prove the left hand side of condition (\ref{eq 1 compactness}), that is,
\begin{equation}\label{Eq : the left hand side of the proposition repetated}
		X \smallsetminus {\uparrow}F^\alpha \sub \bigcup \mathcal{U}^\alpha. 
\end{equation}

Consider $x \in X \smallsetminus {\uparrow}F^\alpha$. We have two cases: either $x \in F^*$ or $x \notin F^*$. First, suppose that $x \in F^*$. By Claim \ref{Claim : cantor compactness}
\[
\text{either} \,\, x \in V_y \text{ for some }y \in G \,\, \text{ or } \, \,  x \in X \smallsetminus {\uparrow}F^\beta \text{ for some }\beta \in \Gamma.
\]
We begin with the case where $x \in V_y$ for some $y \in G$. Since $G \sub \sup^\ast F$ by Claim \ref{Claim : cantor compactness}, we obtain $y \in \sup^\ast F$. Hence, we can apply Claim \ref{Claim : sup Y}(\ref{item : 3 : the new items : comp X}) and the assumption that $x \in V_y$, obtaining 
\[
x \in V_y \subseteq U_y \cup \bigcup \mathcal{W}_y \cup \bigcup \mathcal{U}^{\gamma_y}.
\]
As $y \in G$, the definition of $\mathcal{U}^\alpha$ guarantees that the right hand side of the above display is included in $\bigcup \mathcal{U}^\alpha$. Hence, $x \in \bigcup \mathcal{U}^\alpha$ as desired. Then we consider the case where $x \in X \smallsetminus {\uparrow}F^\beta$ for some $\beta \in \Gamma$. Since $\beta \in \Gamma \subseteq \alpha$, we have $\beta < \alpha$. Therefore, $\beta$ satisfies condition (\ref{eq 1 compactness}). Consequently, from $x \in X \smallsetminus {\uparrow}F^\beta$ it follows that $x \in \bigcup \mathcal{U}^\beta$. As $\beta \in \Gamma$, the definition of $\mathcal{U}^\alpha$ guarantees that $\mathcal{U}^\beta \subseteq \mathcal{U}^\alpha$. Consequently, $\bigcup \mathcal{U}^\beta \subseteq \bigcup \mathcal{U}^\alpha$. Since $x \in \bigcup \mathcal{U}^\beta$, we obtain $x \in \bigcup \mathcal{U}^\alpha$ as desired. This concludes the analysis of the case where $x \in F^*$.

Therefore, we may assume that $x \in X \smallsetminus F^*$. For future reference, it is useful to state the following consequences of this assumption:
\begin{equation}\label{Eq : compactness : one of the final displays}
x \in {\uparrow}F^\beta \text{ for every }\beta < \alpha \, \, \text{ and } \, \, x \notin \sup{\!^{\ast}} F.
\end{equation}

\begin{Claim}\label{Claim : the last claim in the compact proof (hopefully)}
There exist $y^\ast \in \sup^\ast F$ and $z^\ast \in G$ such that $y^\ast \leq x$ and $y^\ast \in V_{z^\ast}$.
\end{Claim}

\begin{proof}[Proof of the Claim]
The left hand side of condition (\ref{Eq : compactness : one of the final displays}) 
guarantees that for each $\beta < \alpha$ there exists $y_\beta \in F^\beta$ such that $y_\beta \leq x$. Since $X$ is a tree and $\alpha$ a limit ordinal, the set $C \coloneqq \{ y_\beta : \beta < \alpha \}$ is a nonempty chain in $F$. Since $C$ is a chain and $X$ a tree, the set
\[
C^* \coloneqq F \cap {\downarrow}C
\]
is also a chain in $F$. Furthermore, from the definition of $y^*$ and $C^*$ it follows that $y^* = \sup C^*$. Therefore, in order to prove that $y^* \in \sup^* F$, it suffices to show that the chain $C^*$ is maximal in $F$. Suppose the contrary, with a view to contradiction. Then there exists $w \in F \smallsetminus C^*$ such that $C^* \cup \{ w \}$ is still a chain. By the definition of $F$ 
there exists $\beta < \alpha$ such that either $w \in F^\beta$ or $w \in X \smallsetminus {\uparrow}\left(F^\beta\right)$. First, suppose that $w \in F^\beta$. Since $w$ and $y_\beta$ are distinct elements of $C^* \cup \{ w \}$, we obtain that either $w < y_\beta$ or $y_\beta < w$. Together with $w, y_\beta \in F^\beta$, this contradicts the assumption that $F^\beta$ is an antichain. Then we consider the case where $w \in X \smallsetminus {\uparrow}F^\beta$. Since $y_\beta \in F^\beta$ and $w$ and $y_\beta$ are two elements of the chain $C^* \cup \{ w \}$, this implies $w < y_\beta$. As $y_\beta \in C$ and $w \in F$, we conclude that $w \in F \cap {\downarrow}C = C^*$, a contradiction. This establishes that the chain $C^*$ is maximal in $F$ and, therefore, $y^* \in \sup^* F$.

It only remains to prove that $y^\ast \in V_{z^\ast}$ for some $z^\ast \in G$. To this end, observe that from $y^\ast \in \sup^\ast F$ and Claim \ref{Claim : cantor compactness} it follows that
\[
y^\ast \in F^* \subseteq \bigcup\limits_{z \in G}
					V_{z}
					\cup
				\bigcup\limits_{\beta \in \Gamma}\left(
					X \smallsetminus {\uparrow} F^{\beta}\right).
\]
Since $y_\beta \in F^\beta$ and $y_\beta \leq y^\ast$ for every $\beta < \alpha$ by construction and $\Gamma \subseteq \alpha$, this yields $y^\ast \in \bigcup_{z \in G}
					V_{z}$. Therefore, there exists $z^\ast \in G$ such that $y^\ast \in V_{z^\ast}$.
\end{proof}

Now, let $y^\ast \in \sup^\ast F$ and $z^\ast \in G \sub \sup^\ast F$ be the elements given by Claim \ref{Claim : the last claim in the compact proof (hopefully)}. Furthermore, let $y_{z^\ast} \in X$ be the element given by Claim \ref{Claim : sup Y}. Lastly, recall that $v_{z^\ast}$ is the element associated with $z^\ast$ at the beginning of the proof of Proposition \ref{proposition : main proposition}.

\begin{Claim}\label{Claim : very last class (hopefully)}
We have that $v_{z^\ast} \leq x$.
\end{Claim}

\begin{proof}[Proof of the Claim]
By Claim \ref{Claim : the last claim in the compact proof (hopefully)} we have $y^\ast \in V_{z^\ast}$. Together with Claim \ref{Claim : sup Y}(\ref{item : 2 : the new items : comp X}), this yields $y^\ast \notin {\uparrow}\left(F^{\gamma_{z^\ast}} \smallsetminus {\uparrow}y_{z^\ast}\right)$. On the other hand, $x \in {\uparrow}F^{\gamma_{z^\ast}}$ by the left hand side of condition (\ref{Eq : compactness : one of the final displays}). Therefore, there exists $w \in F^{\gamma_{z^\ast}}$ such that $w \leq x$. Moreover, $y^\ast \leq x$ by Claim \ref{Claim : the last claim in the compact proof (hopefully)}. Since $X$ is a tree, from $y^\ast, w \leq x$ it follows that $y^\ast$ and $w$ are comparable. Since $y^* \in \sup^\ast F$, the element $y^*$ is the supremum of a maximal chain $C$ in $F$. From the maximality of $C$ and the assumption that $w \in F^{\gamma_{z^\ast}} \subseteq F$, it follows that $y^\ast < w$ is impossible (otherwise $C \cup \{ w \}$ would be a chain in $F$ larger than $C$). Therefore, we conclude that $w \leq y^\ast$. Together with $w \in F^{\gamma_{z^\ast}}$ and $y^\ast \notin {\uparrow}\left(F^{\gamma_{z^\ast}} \smallsetminus {\uparrow}y_{z^\ast}\right)$, this yields $y_{z^\ast} \leq w$. As $w \leq x$, we obtain $y_{z^\ast} \leq x$. Lastly, by Claim \ref{Claim : sup Y}(\ref{item : 1 : the new items : comp X}) we have $v_{z^\ast} \leq y_{z^\ast}$ and, therefore, $v_{z^\ast} \leq x$ as desired.
\end{proof}

We are now ready to conclude the proof, i.e., we establish the left hand side of condition (\ref{Eq : the left hand side of the proposition repetated}) for $x \in X \smallsetminus {\uparrow} F^\alpha$, $x \in X \smallsetminus F^\ast$. We have two cases: either $x \in {\downarrow}Z_{z^\ast}$ or $x \notin {\downarrow}Z_{z^\ast}$. First, suppose that $x \in {\downarrow}Z_{z^\ast}$. Then there exists $w \in Z_{z^\ast}$ such that $x \leq w$. By condition (\ref{Eq : compactness : Vx property}) we have $x \in {\downarrow}w \subseteq \bigcup \mathcal{V}_w$. Since $w \in Z_{z^\ast}$ and $z^\ast \in G$ by 
Claim \ref{Claim : the last claim in the compact proof (hopefully)}, we obtain $\mathcal{V}_w \subseteq \mathcal{U}^\alpha$ and, therefore, $x \in \bigcup \mathcal{U}^\alpha$ as desired. Then we consider the case where $x \notin {\downarrow}Z_{z^\ast}$. Again, we have two cases: either $x \notin {\uparrow}Y_{z^\ast}$ or $x \in {\uparrow}Y_{z^\ast}$. First, suppose that $x \notin {\uparrow}Y_{z^\ast}$. Together with Claim \ref{Claim : very last class (hopefully)}, this yields $x \in {\uparrow}v_{z^\ast} \smallsetminus \left({\uparrow}Y_{z^\ast} \cup {\downarrow}Z_{z^\ast}\right)$. By condition (\ref{Eq : compactness : main lemma property}) this implies $x \in U_{z^\ast}$. As $z^\ast \in G$ by Claim \ref{Claim : the last claim in the compact proof (hopefully)}, the definition of $\mathcal{U}^\alpha$ guarantees that $U_{z^\ast} \in \mathcal{U}^\alpha$. Consequently, $x \in U_{z^\ast} \subseteq \bigcup \mathcal{U}^\alpha$ as desired. 

Lastly, we consider the case where $x \in {\uparrow}Y_{z^\ast}$. We will show that this case never happens, i.e., that it leads to a contradiction. First, there exists $w \in Y_{z^{\ast}}$ such that $w \leq x$. We will prove that $w \in A^\alpha$. Observe that $w \in Y_{z^{\ast}}$ and $z^\ast \in G$ by Claim \ref{Claim : the last claim in the compact proof (hopefully)}. Consequently, to prove that $w \in A^{\alpha}$, it only remains to show that there are no $\beta < \alpha$ and $t \in F^\beta$ such that $w < t$. Suppose, on the contrary, that there exist such $\beta$ and $t$. Recall that $y^\ast \leq x$ by Claim \ref{Claim : the last claim in the compact proof (hopefully)} and that $w \leq x$. Since $X$ is a tree, this yields that $y^\ast$ and $w$ must be comparable. We have two cases: either $y^\ast < w$ or $w \leq y^\ast$. First, suppose that $y^\ast < w$. Together with $w < t$, this yields $y^\ast < t$. Since $y^\ast \in \sup^\ast F$ by Claim \ref{Claim : the last claim in the compact proof (hopefully)}, we know that $y^*$ is the supremum of a maximal chain $C$ in $F$. As $y^* < t $ and $t \in F^{\beta} \subseteq F$, we obtain a contradiction with the maximality of $C$. Then we consider the case where $w \leq y^\ast$. As $w \in Y_{z^\ast}$, we obtain $y^\ast \in {\uparrow}Y_{z^\ast}$. Recall that from Claim \ref{Claim : the last claim in the compact proof (hopefully)} that $y^\ast \in V_{z^\ast}$. Together with $y^\ast \in {\uparrow}Y_{z^\ast}$, this contradicts Claim \ref{Claim : sup Y}(\ref{item : 2 : the new items : comp X}). Hence, we conclude that $w \in A^\alpha$. As the set $A^\alpha$ is finite and $F^\alpha = \min A^\alpha$, from $w \in A^\alpha$ and $w \leq x$ it follows that $x \in {\uparrow}F^\alpha$, contradicting the assumption that $x \in X \smallsetminus {\uparrow}F^\alpha$. This establishes the left hand side of condition (\ref{eq 1 compactness}), thus concluding the argument. 
\end{proof}

\section{Priestley separation axiom}

The aim of this section is to prove the following.

\begin{Theorem}\label{Thm : X is Priestley space : almost done}
The ordered topological space $\langle X; \leq, \tau_{\h{X}} \rangle$ is a Priestley space.
\end{Theorem}

In view of Theorem \ref{Thm : Esakia compactness : trees}, the space $\langle X; \tau_{\h{X}} \rangle$ is compact. Therefore, to establish the above theorem, it suffices to show that $\langle X; \leq, \tau_{\h{X}} \rangle$ satisfies Priestley separation axiom. The rest of this section is devoted to this task.

\begin{Proposition}
The ordered topological space $\langle X; \leq, \tau_{\h{X}} \rangle$ satisfies Priestley separation axiom.
\end{Proposition}

\begin{proof}
We will prove that for every ordinal $\alpha$ and $x, y \in X_{\leq \alpha}$ such that $x \nleq y$ there exists a clopen upset of the ordered topological space $\langle X_{\leq \alpha}; \leq, \tau_{\alpha}\rangle$ such that $x \in U$ and $y \notin U$. The statement will then follow immediately from the case where $\alpha = \h{X}$. During the proof, we will often use $X_{\leq \alpha}$ as a shorthand for $\langle X_{\leq \alpha}; \leq, \tau_{\alpha}\rangle$.
The proof proceeds by induction on $\alpha$.

\subsection*{Base case}

The case where $\alpha = 0$ is straightforward because $X_{\leq \alpha}$ is the singleton  containing the root of $X$.
\subsection*{Successor case}

Suppose that the statement holds for an ordinal $\alpha$ and consider $x, y \in X_{\leq \alpha + 1}$ such that $x \nleq y$. Then for each $z \in \{x,y\} \subseteq X_{\leq \alpha +1}$ let
	\[
		\bar{z} \coloneqq 
			\begin{cases}
				z & \text{ if } z \in X_{\leq \alpha};\\
				\text{the immediate predecessor of }z & \text{ if }z \in X_{\alpha + 1}.
			\end{cases}		
	\]
	Clearly, $\bar{z} \leq z$ and $z \in X_{\leq \alpha}$. We have two cases: either $\bar{x} \nleq \bar{y}$ or $\bar{x} \leq \bar{y}$.

First, suppose that $\bar{x} \nleq \bar{y}$. Since $\bar{x}, \bar{y} \in X_{\leq \alpha}$ and $\bar{x} \nleq \bar{y}$, we can apply the inductive hypothesis obtaining a clopen upset $V$ of $X_{\leq \alpha}$ such that $\bar{x} \in V$ and $\bar{y} \notin V$.  Then let
	\[
		U \coloneqq V \cup {\uparrow_{\alpha +1}}\left(V \cap X_\alpha\right).
	\]

We will prove that $U$ is a clopen upset of $X_{\leq \alpha +1}$. From the assumption that $V$ is an upset of $X_{\leq \alpha}$ it follows that $U$ is an upset of $X_{\leq \alpha+1}$. Furthermore, the fact that $V$ is an open set of $X_{\leq \alpha}$ and the definition of $\mathcal{S}_{\alpha +1}$ guarantee that $U$ is an open set of $X_{\leq \alpha +1}$. Therefore, it only remains to show that $U$ is a closed set of $X_{\leq \alpha +1}$. Since $V$ is an upset of $X_{\leq \alpha}$ and $X$ a tree, we have
\[
		X
			\smallsetminus
				\left(
					V \cup {\uparrow}\left(V \cap X_{\alpha}\right) 
				\right)
			=
		\left(X_{\leq \alpha} \smallsetminus V\right)
					\cup
		{\uparrow}
						\left(X_{\alpha} \smallsetminus V\right). 
	\]
Using the definition of $U$ and restricting to $X_{\leq \alpha +1}$ both sides of the above equality, we obtain
\[
	X_{\leq \alpha +1} \smallsetminus U =	X_{\leq \alpha +1}
			\smallsetminus
				\left(
					V \cup {\uparrow}_{\alpha +1}\left(V \cap X_{\alpha}\right) 
				\right)
			=
\left(X_{\leq \alpha} \smallsetminus V\right)
					\cup
				{\uparrow_{\alpha+1}}
						\left(\left(X_{\leq \alpha} \smallsetminus V\right) \cap X_{\alpha}\right).
	\]
As $X_{\leq \alpha} \smallsetminus V \in \tau_\alpha$ by assumption, the definition of $\mathcal{S}_{\alpha+1}$ guarantees that the right hand side of the above display is an open set of $X_{\leq \alpha +1}$. Hence, $U$ is a closed set of $X_{\leq \alpha +1}$. This establishes that $U$ is a clopen upset of $X_{\leq \alpha +1}$.

Therefore, it only remains to prove that $x \in U$ and $y \notin U$. Recall that $\bar{x} \in V$  and $\bar{x} \leq x \in X_{\leq \alpha +1}$.\ As $U$ is an upset of $X_{\leq \alpha +1}$ containing $V$, we obtain $x \in U$.  To prove that $y \notin U$, we consider separately two cases: either $y \in X_{\leq \alpha}$ or $y \in X_{\alpha +1}$. First suppose that $y \in X_{\leq \alpha}$. Then $y = \bar{y}$. As $\bar{y} \notin V$ by assumption, we also have $y \notin V$. Together with $y \in X_{\leq \alpha}$, this yields $y \notin V \cup {\uparrow}_{\alpha +1} \left(V \cap X_\alpha\right) = U$ as desired. Then we consider the case where $y \in X_{\alpha +1}$. We have $y \notin V$ because $V \subseteq X_{\leq \alpha}$. Moreover, $y \notin {\uparrow}_{\alpha+1}\left(V \cap X_\alpha\right)$ because $\bar{y}$, which is the only predecessor of $y$ of height $\alpha$, does not belong to $V$ by assumption. Hence, we conclude that $y \notin U$.

It only remains to consider the case where $\bar{x} \leq \bar{y}$. As $\bar{y} \leq y$ and $x \nleq y$, we have $\bar{x} \neq x$. By the definition of $\bar{x}$ this implies $x \in X_{\alpha+1}$ and $\bar{x} \in X_\alpha$. Therefore, from $\bar{x} \leq \bar{y} \in X_{\leq \alpha}$ it follows that $\bar{x} = \bar{y}$. Hence, $\bar{y}  = \bar{x} \in X_{\alpha}$.  We have three subcases: 
\[
\text{either }\, \, y \in P_\alpha\, \, \text{ or }\, \, y \in S_{\alpha+1}\, \, \text{ or }\, \, y \notin P_\alpha \cup S_{\alpha+1}.
\]

Suppose first that $y \in P_\alpha$. We will prove that ${\downarrow}y$ is a clopen set of $X_{\leq \alpha +1}$. Since $y \in P_\alpha$, the definition of $\mathcal{S}_{\alpha +1}$ guarantees that ${\downarrow}y$ is an open set of $X_{\leq \alpha +1}$. By the same token the set
\[
\left(X_{\leq \alpha} \cup {\uparrow_{\alpha+1}}\left(X_{\leq \alpha}\cap X_\alpha\right)\right) \smallsetminus {\downarrow}y
\]
is also an open of $X_{\leq \alpha +1}$, which is easily seen to coincide with $X_{\leq \alpha +1} \smallsetminus {\downarrow}y$. Therefore, ${\downarrow}y$ is a clopen set of $X_{\leq \alpha +1}$. Together with $x \nleq y$, this implies that $X_{\leq \alpha+1} \smallsetminus {\downarrow}y$ is a clopen upset of $X_{\leq \alpha+1}$ containing $x$ but not $y$ and we are done.

Then we consider the case where $y \in S_{\alpha+1}$. As before, it suffices to show that ${\downarrow}y$ is a clopen set of $X_{\leq \alpha +1}$. The fact that ${\downarrow}y$ is closed is proved as in the previous case. To prove that it is open, observe that $\bar{y} \in P_\alpha$ because $y \in S_{\alpha+1}$. From the definition of $\mathcal{S}_{\alpha+1}$ and the assumption that $y \in S_{\alpha+1}$ and $\bar{y} \in P_\alpha$ it follows that both $\{ y \}$ and ${\downarrow}\bar{y}$ are open sets of $X_{\leq \alpha +1}$. Therefore, ${\downarrow}y = \{ y \} \cup {\downarrow}\bar{y}$ is an open set of $X_{\leq \alpha +1}$ as desired.

Lastly, we consider the case where $y \notin P_\alpha \cup S_{\alpha+1}$.\ We will prove that $x \in S_{\alpha+1}$. Suppose the contrary, with a view to contradiction. As $x \in X_{\alpha+1}$ and $\bar{x} \in X_\alpha$, from $x \notin S_{\alpha+1}$ it follows that $x = \bar{x}^+$. Moreover, from $\bar{y} = \bar{x} \in X_\alpha$ and $\bar{y} \leq y \in X_{\leq \alpha+1}$ it follows that either $y \in \{ \bar{y}, \bar{y}^+ \} \cup S_{\alpha+1}$. As $y \notin S_{\alpha+1}$ by assumption, we get $y \in \{ \bar{y}, \bar{y}^+ \}$. Moreover, from $\bar{y} = \bar{x} < x$ and $\bar{y} \in X_\alpha$ it follows that $\bar{y} \in P_\alpha$. Together with $y \in \{ \bar{y}, \bar{y}^+ \}$ and the assumption that $y \notin P_\alpha$, this yields $y = \bar{y}^+$. Since $\bar{x} = \bar{y}$ and $\bar{x}^+ = x$, we obtain $y = x$, a contradiction with $x \nleq y$. Hence, we conclude that $x \in S_{\alpha +1}$.

We will use this fact to prove that $\{ x \}$ is a clopen upset of $X_{\leq \alpha+1}$ containing $x$ but not $y$. As $x \nleq y$ and $x$ is a maximal element of $X_{\leq \alpha+1}$ (the latter because $x \in X_{\alpha+1}$), it suffices to show that $\{ x \}$ is a clopen set of $X_{\leq \alpha +1}$. Since $x \in S_{\alpha+1}$, the definition of $\mathcal{S}_{\alpha+1}$ guarantees that $\{ x \}$ is an open set of $X_{\leq \alpha +1}$. To prove that it is also closed, observe that
	\[
		X_{\leq \alpha +1} \smallsetminus \{x\}
			=
				\left(
					\left(
						X_{\leq \alpha}
							\cup
						{\uparrow_{\alpha+1}}
							\left(
								X_{\leq \alpha} \cap X_\alpha
							\right)
					\right)
						\smallsetminus {\downarrow} x 
				\right)
				\cup {\downarrow}\bar{x}
	\]
because $x \in X_{\alpha+1}$ and $\bar{x}$ is the unique immediate predecessor of $x$. Furthermore, as $\bar{x} \in P_\alpha$ and $x \in S_{\alpha+1}$, the right hand side of the above display is the union of two members of $\mathcal{S}_{\alpha+1}$. Hence, $\{ x \}$ is a closed set of $X_{\leq \alpha+1}$ as desired.

\subsection*{Limit case}

Suppose that $\alpha$ is a limit ordinal and consider $x, y \in X_{\leq \alpha}$ such that $x \nleq y$. We will prove that there exist $\beta < \alpha$ and $x^\ast \in X_{\leq \beta}$ such that $x^\ast \leq x$ and $x^\ast \nleq y$. If $x \in X_{< \alpha}$, we are done letting $x^\ast \coloneqq x$ and $\beta \coloneqq \h{x}$. Then we consider the case where $x \in X_\alpha$. Since $\alpha$ is a limit ordinal and every nonempty chain in $X$ has a supremum, from $x \in X_\alpha$ it follows that $x$ is the supremum of the nonempty chain ${\downarrow}x \smallsetminus \{ x \}$. As $x \nleq y$, this implies that there exists $x^\ast \in {\downarrow}x \smallsetminus \{ x \}$ such that $x^\ast \nleq y$. Letting $\beta \coloneqq \h{x^*}$ and observing that $\beta < \alpha$, we are done. 

Now, consider the nonempty chain $C \coloneqq X_{\leq \beta} \cap {\downarrow}y$. By assumption the supremum $y^*$ of $C$ exists and, moreover, belongs to $X_{\leq \beta}$ because $C \subseteq X_{\leq \beta}$. Since $x^\ast \nleq y$ and $y^* \leq y$, we have $x^\ast \nleq y^\ast$. Recall that $\beta < \alpha$. As $x^\ast, y^\ast \in X_{\leq \beta}$ and $x^\ast \nleq y^\ast$, the inductive hypothesis guarantees the existence of a clopen upset $U$ of $X_{\leq \beta}$ such that $x^\ast \in U$ and $y^\ast \notin U$. Since $\alpha$ is a limit ordinal, the definition of $\mathcal{S}_\alpha$ ensures that both
\[
U \cup {\uparrow}_\alpha\left(U \cap X_\beta\right) \, \, \text{ and } \, \, \left(X_{\leq \beta} \smallsetminus U\right) \cup {\uparrow}_\alpha\left(X_\beta \smallsetminus U\right)
\]
are open sets of $X_{\leq \alpha}$. As $U$ is an upset of $X_{\leq \beta}$, the set on the left hand side of the above display coincides with ${\uparrow}_\alpha U$. Similarly, the set of the right hand side of the display is $X_{\leq \alpha} \smallsetminus {\uparrow}_\alpha U$ because $X$ is a tree and $U$ an upset of $X_{\leq \beta}$. Therefore, ${\uparrow}_\alpha U$ is a clopen upset of $X_{\leq \alpha}$. 

Lastly, from $x^\ast \in U$ and $x^\ast \leq x \in X_{\leq \alpha}$ it follows that $x \in {\uparrow}_\alpha U$. Therefore, it only remains to prove that $y \notin {\uparrow}_\alpha U$. Since ${\uparrow}_\alpha U = U \cup {\uparrow}_\alpha\left(U \cap X_\beta\right)$, it suffices to show that $y \notin U$ and $y \notin {\uparrow}_\alpha\left(U \cap X_\beta\right)$. 
Suppose the contrary, with a view to contradiction. We have two cases: either $y \in U$ or $y \in {\uparrow}_\alpha\left(U \cap X_\beta\right)$. First, suppose that $y \in U$. Then $y = y^\ast$ because $y \in U \subseteq X_{\leq \beta}$ and $y^\ast$ is the supremum of ${\downarrow} y \cap X_{\leq \beta}$.\ But this implies $y^\ast = y \in U$, which is false. Then we consider the case where $y \in {\uparrow}_\alpha\left(U \cap X_\beta\right)$. The definition of $y^\ast$ guarantees that $y^\ast \in U \cap X_\beta$, a contradiction with $y^\ast \notin U$. Hence, we conclude that $y \notin {\uparrow}_\alpha U$.
\end{proof}

\section{The end}
In order to conclude the proof of Theorem \ref{Thm:main-forest}, we need to show that $\langle X; \leq, \tau_{\h{X}} \rangle$ is an Esakia space. As $\langle X; \leq, \tau_{\h{X}} \rangle$ is a Priestley space by Theorem \ref{Thm : X is Priestley space : almost done}, it only remains to prove that the downset of every open set is still open. Therefore, the following observation concludes the proof of Theorem \ref{Thm:main-forest}.

\begin{Proposition}
	For every $U \in \tau_{\h{X}}$ we have ${\downarrow}U \in \tau_{\h{X}}$.
\end{Proposition}

\begin{proof} 
The proof hinges on the following claim:

\begin{Claim}\label{lemma : nonmaximal singletons}
	Let $\alpha$ be an ordinal and $x \in X_{\leq\alpha} \smallsetminus \max X_{\leq \alpha}$. Then ${\downarrow} x \in \mathcal{S}_{\alpha}$. 
\end{Claim}

\begin{proof}[Proof of the Claim]
The proof of the claim proceeds by induction on $\alpha$.

\subsection*{Base case} The case where $\alpha = 0$ holds vacuously because $X_{\leq 0} \smallsetminus \max X_{\leq 0} = \emptyset$ and, therefore, $x \in X_{\leq 0} \smallsetminus \max X_{\leq 0}$ is impossible.

\subsection*{Successor case} Suppose that $x \in X_{\leq\alpha+1} \smallsetminus \max X_{\leq \alpha+1}$. We have two cases: either $x \in \max X_{\leq \alpha}$ or $x \notin \max X_{\leq \alpha}$. First, suppose that $x \in \max X_{\leq \alpha}$. Since $x \notin \max X_{\leq \alpha+1}$, this implies $x \in P_\alpha$. Consequently, ${\downarrow}x \in \mathcal{S}_{\alpha+1}$ by the definition of $\mathcal{S}_{\alpha+1}$.\ Then we consider the case where $x \notin \max X_{\leq \alpha}$. Together with the assumption that $x \in X_{\leq\alpha+1} \smallsetminus \max X_{\leq \alpha+1}$, this yields $x \in X_{< \alpha}$. As $x \notin \max X_{\leq \alpha}$, we can infer $x \in X_{\leq \alpha} \smallsetminus\max X_{\leq \alpha}$. Consequently, we can apply the inductive hypothesis, obtaining ${\downarrow}x \in \mathcal{S}_\alpha$. By the definition of $\mathcal{S}_{\alpha+1}$ we have
\[
{\downarrow}x \cup {\uparrow}_{\alpha+1}(X_{\alpha} \cap {\downarrow}x) \in \mathcal{S}_{\alpha+1}.
\]
Furthermore, from $x \in X_{< \alpha}$ it follows that $X_{\alpha} \cap {\downarrow}x = \emptyset$. Therefore, the above display simplifies to ${\downarrow}x \in \mathcal{S}_{\alpha+1}$ and we are done.

\subsection*{Limit case} Let $x \in X_{\leq\alpha} \smallsetminus \max X_{\leq \alpha}$ and assume that $\alpha$ is a limit ordinal. As $x \notin \max X_{\leq \alpha}$, we have $\h{x} < \alpha$. We will prove that
\[
x \in X_{\leq \h{x}+1} \smallsetminus \max X_{\leq \h{x}+1}.
\]
It is clear that $x \in X_{\leq \h{x}+1}$. Therefore, it suffices to prove that $x \notin \max X_{\leq  \h{x}+1}$. Suppose the contrary, with a view to contradiction. From $x \in \max X_{\leq \h{x}+1}$ and the fact that $x$ has order type $\h{x}$ it follows that $x$ is a maximal element of $X$. Together with $\h{x} \leq \alpha$, this yields $x \in \max X_{\leq \alpha}$, a contradiction. This establishes the above display.

Recall that $\h{x} < \alpha$. Since $\alpha$ is a limit ordinal, this yields $\h{x} +1 < \alpha$. Therefore, we can apply the inductive hypothesis to the above display, obtaining ${\downarrow}x \in \mathcal{S}_{\h{x} +1}$. Since $\alpha$ is a limit ordinal, the definition of $\mathcal{S}_\alpha$ guarantees that 
\[
{\downarrow}x \cup {\uparrow}_\alpha(X_{\h{x}+1} \cap {\downarrow}x) \in \mathcal{S}_\alpha.
\]
As ${\downarrow}x \subseteq X_{\leq \h{x}}$, we have $X_{\h{x}+1} \cap {\downarrow}x = \emptyset$. Therefore, the above display simplifies to ${\downarrow}x \in \mathcal{S}_\alpha$.
\end{proof}

Now, we turn to prove the main statement. Let $U \in \tau_{\h{X}}$. Clearly, we have
	\[
		{\downarrow} U
			=
		U
			\cup 
		\bigcup \{ {\downarrow}w : w \in {\downarrow}U \colon w \notin \max X\}.
	\]
	As $U \in \tau_{\h{X}}$ by assumption and ${\downarrow}w \in \tau_{\h{X}}$ for each $w \notin \max X$ by  Claim \ref{lemma : nonmaximal singletons}, the right hand side of the above display belongs to the topology $\tau_{\h{X}}$. Hence, we conclude that ${\downarrow}U \in \tau_{\h{X}}$.
\end{proof}

\vspace{.5cm}
\paragraph{\bfseries Acknowledgements.}
The first author was supported by the FPI scholarship PRE$2020$-$093696$, associated with the I+D+i research project PID$2019$-$110843$GA-100 I+D \textit{La geometría de las lógicas no clásicas} funded by the Ministry of Science and Innovation of Spain.\ The second author was supported by the proyecto PID$2022$- $141529$NB-C$21$ de investigaci\'on financiado por MICIU/AEI/$10$.$13039$/$501100011033$ y por FEDER, UE. He was also supported by the Research Group in Mathematical Logic, $2021$SGR$00348$ funded by the Agency for Management of University and Research Grants of the Government of Catalonia. All the authors were supported by the MSCA-RISE-Marie Skłodowska-Curie Research and Innovation Staff Exchange (RISE) project MOSAIC $101007627$ funded by Horizon $2020$ of the European Union.

\bibliographystyle{plain}

\end{document}